\documentclass[12pt]{article}
\usepackage[final]{epsfig}
\usepackage{amsfonts, amsmath, amssymb, amsthm}
\usepackage{color}
\usepackage{graphicx}

\numberwithin{equation}{section}

\newtheorem{theorem}{Theorem}[section]
\newtheorem{lemma}{Lemma}[section]

\newcommand{\R}{\mathbb{R}} 
\newcommand{\C}{\mathbb{C}}

\newcommand{\N}{\mathbb{N}}

\newcommand{\be}{\begin{equation}}
\newcommand{\ee}{\end{equation}}
\newcommand{\calA}{{\mathcal{A}}}
\newcommand{\calB}{{\mathcal{B}}}
\newcommand{\calM}{{\mathcal{M}}}
\newcommand{\calE}{{\mathcal{E}}}
\newcommand{\eps}{{\varepsilon}}

\newcommand{\weakly}{{\rightharpoonup}}
\newcommand{\phinu}{{\phi^{(\nu)}}}
\newcommand{\munu}{{\mu^{(\nu)}}}
\newcommand{\phinuH}{{\phi_H^{(\nu)}}}
\newcommand{\phinuL}{{\phi_L^{(\nu)}}}
\newcommand{\phinui}{{\phi_i^{(\nu)}}}

\begin{document}

\begin{center}

\begin{Large}
   New optimal control problems in density functional theory  \\
      motivated by photovoltaics
\end{Large}
\normalsize
\vspace{5mm}

Gero Friesecke$^1$ and Michael Kniely$^2$ \\[1mm]
$^1\,$Faculty of Mathematics, Technische Universit\"at M\"unchen, {\tt gf@ma.tum.de}  \\
$^2\,$Institute of Science and Technology Austria, {\tt michael.kniely@ist.ac.at} \\[2mm]
 
\end{center}


                                                                                                      
{\bf Abstract.} We present and study novel optimal control problems motivated by the search for photovoltaic materials with high power-conversion efficiency. The material must perform the first step: convert light (photons) into electronic excitations. We formulate various desirable properties of the excitations as mathematical control goals at the Kohn-Sham-DFT level of theory, with the control being given by the nuclear charge distribution. We prove that nuclear distributions exist which give rise to optimal HOMO-LUMO excitations, and present illustrative numerical simulations  
for 1D finite nanocrystals. We observe pronounced goal-dependent features such as large electron-hole separation, and a hierarchy of length scales: internal HOMO and LUMO wavelengths $<$ atomic spacings $<$ (irregular) fluctuations of the doping profiles $<$ system size. 
\section{Introduction} \label{sec:Intro} 
\begin{figure} 
\begin{center}
\includegraphics[width=\textwidth, trim=50 130 35 120, clip]{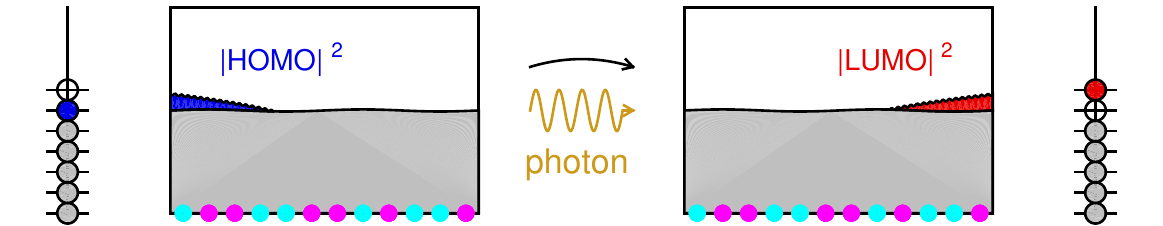}
\end{center}
\caption{Schematic picture of 
light-induced electron-hole pair formation.} 
\label{F:schematic}
\end{figure}

In this paper we propose and study novel optimal control problems motivated by the ongoing search for photovoltaic materials with high power-conversion efficiency (see, e.g., \cite{JN12, KOBH13}). 
Photovoltaic devices convert light (photons) into current (electron motion). The photovoltaic material must perform the first step: to convert light (photons) into suitable electronic excitations. 
\\[2mm]
In heterojunction solar cells, the photovoltaic material is placed inside a heterojunction 
and the electronic excitations of interest are typically electron-hole-pairs forming inside the junction. These excitations must have a highly nontrivial list of properties: long lifetime; appropriate size of energy gap as compared to the electronic ground state; spatial electron-hole-separation; high mobility. These properties depend critically on the type of atoms in the material (e.g., in a crystalline material, the doping profile, or, in a polymer, the heteroatom substitutions \cite{LK88}); see \cite{KOBH13} for a systematic account. 
\\[2mm]
The typical procedure in the ab-initio modelling of materials, including photovoltaic ones \cite{ZLL11, MW14, HLYDY16}, is to begin from the real material (specified, e.g., by the constituent atoms and their spatial arrangement), select a suitable mathematical model for the relevant quantities (e.g., electronic excitations), and 
simulate the model to predict behaviour. Here we take first steps towards the {\it inverse} route: begin with a desirable property or effect; formulate it as a mathematical property of solutions to the governing equations; then try to find a set of parameters (nuclear charges and positions) which produces this effect. In our context of photovoltaic materials, a suitable model for the electronic structure is given by the Kohn-Sham equations (see \eqref{KS}--\eqref{KSeqs} below) and 
electron-hole pair excitations can, as a first approximation, be taken to correspond to HOMO-LUMO transitions (see \eqref{HL} and Figure \ref{F:schematic}), and we are thus led 
to novel optimal control problems for the KS equations (see Section \ref{S:OCPs}) which appear to be very interesting both from a physical and a mathematical point of view. 
Our goal in this paper is threefold: 
\begin{itemize}
\item[--] 
to introduce these control problems and in particular formulate various desirable properties of photovoltaic materials as mathematical control goals, with the control being given by ``the choice of material'', i.e., the nuclear charge distribution (see Section 2)
\item[--] 
to place the problems on a firm mathematical footing and prove rigorously that they are well posed, i.e. that optimal nuclear charge distributions and ensuing optimal electronic excitations exist (see Section 3)
\item[--]
to present illustrative numerical simulations of optimal doping profiles and the resulting HOMO-LUMO-excitations for one-dimensional finite nanocrystals with 20 atoms and 120 electrons (numbers chosen so that the un-doped chain corresponds to pure carbon).
\end{itemize}
The numerical optimizers, described in detail in Section 4, are seen to exhibit a remarkable and nontrivial multi-scale structure, with (in atomic units of length, 1 a.u. $\approx 0.5 \cdot 10^{-10}$ m) \\
-- atomic spacings $\sim$ 1 \\
-- internal HOMO and LUMO wavelengths $\sim 0.2$ \\
-- doping profiles which are irregular and have wavelength $\sim 2$ or $3$ \\
-- diameter and spatial separation of HOMO and LUMO $\sim 10$ \\ 
-- total system diameter $\sim 20$. 
\section{Desirable properties of opto-electronic excitations as mathematical control goals}
We begin by briefly recalling Kohn-Sham density functional theory \cite{HK64, KS65, PY89}, then descibe how to model excitations at the KS-DFT level of theory, then formulate various optimal control problems associated with excitations.
\subsection{Kohn-Sham equations} \label{S:KS}
The standard electronic structure model in materials science is Kohn-Sham density functional theory (KS-DFT), due to its good compromise between accuracy and computational feasibility for large systems. In our context of photovoltaics we assume that the electrons are located in an open bounded region $\Omega\subset\R^3$ (the photovoltaic element) and the atomic nuclei are distributed in a compact subset $\Omega_{nuc}\subset\Omega$. The nuclear charge distribution will be denoted by $\mu$; the prototypical example is a finite number of point charges with charges $Z_\alpha>0$ and positions $R_\alpha\in\Omega_{nuc}$, i.e. $\mu = \sum_{\alpha=1}^M Z_\alpha\delta_{R_\alpha}$. For simplicity we make the customary assumptions that the number of electrons is even, i.e. equal to $2n$ for some $n\in\N$, and that the system is not spin-polarized. In this case, KS-DFT models the electrons by $n$ orbitals $\phi_1,...,\phi_n\, : \, \R^3\to\C$ 
(the physical picture is that each orbital is occupied by two electrons of opposite spin). The orbitals must be $L^2$-orthonormal, i.e.
\begin{equation} \label{ortho}
      \langle \phi_i, \, \phi_j \rangle = \delta_{ij} \;\;\;\;\;\; (i, \, j = 1,...,n)
\end{equation}
where $\langle \phi_i, \, \phi_j\rangle = \int_\Omega \overline{\phi_i(x)} \, \phi_j(x) \, dx$ denotes the $L^2$ inner product, and in the ground state the collection $\phi=(\phi_1,..,\phi_n)$ of orbitals is governed by the Kohn-Sham variational principle:
\begin{equation} \label{KS}
  \phi \in \mbox{argmin}\, \calE_\mu \; \mbox{ subject to the constraints \eqref{ortho}},  
\end{equation}
where 
\begin{equation} \label{EKS}
  \calE_\mu[\phi_1,...,\phi_n] =  \underbrace{\sum_{i=1}^n 2\int_\Omega\! \frac12|\nabla\phi_i|^2}_{=:T[\phi]}  + \underbrace{\int_\Omega\! v_{ext}\rho}_{=:V_{ext}[\rho ]} + \underbrace{\frac12\int\int_{\Omega\times\Omega}\!\! 
    \frac{\rho(x)\rho(y)}{|x-y|}dx\, dy}_{=:J_H[\rho]} + \underbrace{\int_\Omega\! e_{xc}(\rho(x)) \, dx}_{=:E_{xc}[\rho]}
\end{equation}
is the Kohn-Sham energy functional and
\be \label{density}
      \rho(x)=2 \sum_{i=1}^n |\phi_i(x) |^2 
\ee 
is the total electron density. Minimizers satisfy the Kohn-Sham equations 
\be \label{KSeqs}
    \underbrace{\Bigl(-\frac12 \Delta + v_{ext} + v_H + v_{xc}(\rho)\Bigr)}_{=:h_\phi} \phi_i = \sum_{j=1}^n \lambda_{ij}\phi_j
\ee
where the $\lambda_{ij}$ are Lagrange multipliers coming from the orthonormality constraints \eqref{ortho}. The meaning of the different terms in the energy functional and the KS equations is the following: $T$ is the kinetic energy; $V_{ext}$ is the electron-nuclei energy with $v_{ext}$ being the electrostatic potential of the nuclei, see \eqref{pots}; 
$J_H$ (the Hartree energy) would correspond to the interelectron repulsion energy if the electrons were mutually independent; $E_{xc}$ (the exchange-correlation energy) is a correction accounting for correlation effects. The function $e_{xc}$ is a ``known'' function of $\rho$ (its pointwise values 
$e_{xc}(\overline{\rho})$ ($\overline{\rho}\ge 0$) model the exchange-correlation energy per unit volume of a homogeneous electron gas with density $\overline{\rho}$). Modelling the exchange-correlation energy in this way corresponds to the local density approximation (LDA). A simple prototype which accounts only for exchange is the Dirac exchange energy 
\be \label{dirac}
      e_{xc}(\rho)= - c_x \rho^{4/3}, \;\;\; c_x = \tfrac{3}{4} \Bigl(\tfrac{3}{\pi}\Bigr)^{\frac{1}{3}}. 
\ee
The potentials in \eqref{KSeqs} are the external, Hartree, and exchange-correlation potentials,
\be \label{pots}
   v_{ext}(x) = \! -\! \int_{\Omega_{nuc}} \! \frac{1}{|x-y|} d\mu(y), \;\;\;
   v_H(x)\! =\! \int_{\Omega}\! \frac{\rho(y)}{|x-y|} dy, \;\;\;
   v_{xc}(\rho) \!=\! \tfrac{d}{d\rho} e_{xc}(\rho),
\ee
and $\mu$ is the nuclear charge density, assumed to satisfy 
\be \label{neutrality}
    \int_{\Omega_{nuc}} d\mu =2n \;\; \mbox{(charge neutrality)}.
\ee 
The effective one-body operator $h_\phi$ in \eqref{KSeqs} (the Kohn-Sham Hamiltonian) depends also on the nuclear charge density $\mu$ (which is suppressed by our notation). This operator is invariant under unitary transformations $\phi_i \mapsto \phi_i'=\sum_j U_{ij}\phi_j$, $U$ a unitary $n\times n$ matrix, and by a suitable unitary transformation of the orbitals the KS equations \eqref{KSeqs} can be brought into the canonical form
\be \label{cKSeqs}
   h_\phi \varphi_i = \eps_i \varphi_i \;\;\; (\eps_1\le \eps_2 \le ... \le \eps_n ).
\ee
For common choices of the exchange-correlation functional the $\eps_i$ correspond to the lowest $n$ eigenvalues of $h_{\phi}$, accounting for multiplicity (for a rigorous proof of this fact see \cite{FG18}). 
\subsection{Excitations}
{\bf HOMO-LUMO-transition} We limit ourselves here to the lowest opto-electronic excitation of the system, and the most basic model for it within KS-DFT, the HOMO-LUMO-transition. As we will see, the resulting control problems are already very interesting and highly nontrivial both physically and mathematically. Treatment of the whole excitation spectrum, as well as of many-body corrections like the Casida ansatz, lies beyond the scope of this paper. 

In the HOMO-LUMO transition, 
an electron pair migrates from the highest occupied molecular orbital (HOMO) to the lowest unoccupied molecular orbital (LUMO),  
\be \label{HL}
   (\phi_1,...,\phi_{n-1},\textcolor{blue}{\phi_n}) \;\; \longrightarrow \;\; (\phi_1,...,\phi_{n-1},\textcolor{red}{\phi_{n+1}}).
\ee
The transition is induced by an incoming photon whose frequency $\nu$ satisfies the Bohr condition $h \nu = \eps_{n+1}-\eps_n$, where $h$ is Planck's constant. 
Here the orbitals $\phi_1,...,\phi_n$ are the canonical KS orbitals (i.e., the eigenstates in \eqref{cKSeqs} ordered by size of eigenvalue), $\phi_n$ is the HOMO, and $\phi_{n+1}$ (the LUMO) is the next
eigenstate of $h_\phi$ (i.e. the eigenstate $L^2$-orthogonal to $\phi_1,...,\phi_n$ with lowest eigenvalue). According to a recent study \cite{P17}, the HOMO-LUMO eigenvalue difference $\eps_{n+1}-\eps_n$ correctly reproduces the KS electron-hole pair creation energy (i.e. the difference between ionization energy and electron affinity) to O($\tfrac{1}{N}$) if the density change by adding an electron or hole is delocalized, where $N$ is the number of electrons in the system. 

Of course, it may happen (even though we did not observe it numerically for {\it optimal} excitations) that HOMO and/or LUMO are nonunique, due to eigenvalue crossings, and we allow for this in our analysis. 
\\[2mm]
{\bf Variational definition of HOMO and LUMO} 
The following variational definition works irrespective of degeneracies, and will be very convenient for the mathematical analysis of optimal excitations. Consider the quadratic form associated with the KS Hamiltonian $h_\phi$ from \eqref{KSeqs}, 
\begin{equation} \label{EHL}
   \calE_{\mu,\phi}[\chi] = \langle \chi, \, h_\phi \chi\rangle = 
                  \frac12 \int_\Omega |\nabla\chi|^2 + \int_\Omega \Bigl(v_{ext} + v_H + v_{xc}(\rho) \Bigr) \rho_\chi, \;\;\; \rho_\chi = |\chi|^2,         
\end{equation}
where $v_{ext}$, $v_H$, $v_{xc}(\rho)$ are the usual potentials recalled in \eqref{pots}. We call this quadratic single-particle functional the {\it excitation functional}. The potential $v_{ext}$ depends on the nuclear charge density $\mu$, and the other two potentials  depend on the single-particle density \eqref{density} associated with the occupied KS orbitals $\phi=(\phi_1,...,\phi_n)$. Thus the excitation functional depends on both the nuclear charge density $\mu$ and the occupied KS orbitals $\phi$, as emphasized by the notation $\calE_{\mu,\phi}$. We now define a HOMO $\varphi_H$ by
\begin{equation} \label{VHOMO}
    \varphi_H \in \mbox{argmax}\, \calE_{\mu,\phi}   \mbox{ subject to the constraints }
    \varphi_H\in \mbox{Span}\, \{ \phi_1,...,\phi_n \}, \;\, \langle \varphi_H,\, \varphi_H\rangle = 1 
\end{equation} 
and a LUMO $\varphi_L$ by
\begin{equation} \label{VLUMO}
    \varphi_L \in \mbox{argmin}\, \calE_{\mu,\phi}   \mbox{ subject to the constraints }
    \langle\phi_i, \, \varphi_L\rangle = 0 \; (i=1,...,n), \;\, \langle \varphi_L,\, \varphi_L\rangle = 1.
\end{equation}    
Solutions to these variational problems -- i.e., a HOMO and a LUMO -- can be proven to exist (see Section 3), and obviously satisfy the KS equations
\be \label{HOMOLUMOeqs}
   h_\phi \varphi_H = \eps_H \varphi_H, \;\;\; h_\phi \varphi_L = \eps_{L} \varphi_L
\ee
for some eigenvalues $\eps_H$ (the HOMO energy) and $\eps_L$ (the LUMO energy). As discussed below eq.~\eqref{cKSeqs}, for common choices of the exchange-correlation functional we have $\eps_H=\eps_n$, $\eps_L=\eps_{n+1}$, where $\eps_1\le \eps_2 \le ... $ are the eigenvalues of $h_\phi$.  
\subsection{Optimal control problems} \label{S:OCPs}
We limit ourselves here to four basic optimal control problems. 
Of course, many variants and especially combinations of these can be considered. These control problems are associated with different physical properties of excitations which are desirable in photovoltaics: \\
-- charge transfer over a large distance \\
-- no spatial electron--hole overlap \\
-- long lifetime \\
-- prescribed bandgap. 

Mathematically, the ensuing optimal control problems for excitations correspond to certain novel {\it optimal control problems with PDE constraints}: 

\begin{center}
\begin{minipage}{0.85\textwidth}
Choose the nuclear charge density $\mu$ (the {\it control}) so as to minimize a functional
\be \label{goal}
    J [\phi, \, \phi_H, \, \phi_L, \mu]
\ee
(the {\it control goal}) subject to the fact that $\phi$ must be the Kohn-Sham ground state for the nuclear charge density $\mu$ and $\phi_H$, $\phi_L$ are the associated HOMO and LUMO, i.e.
\begin{equation} \label{stateequation}
    \phi \mbox{ satisfies \eqref{KS}}, \;\;\; \phi_H \mbox{ satisfies \eqref{VHOMO}}, \;\;\; 
    \phi_L \mbox{ satisfies \eqref{VLUMO}}
\end{equation}
(the {\it state equation}). 
\end{minipage}
\end{center}

Different desirable physical properties of opto-electronic excitations are to be modelled by different mathematical control goals $J$; see our list of examples below. Note that the control, $\mu$, often does not appear explicitly in the control goal, but enters in a highly indirect and nonlinear way through the state equation.
\\[2mm] 
{\bf a) Charge transfer.} The following functional measures the amount of charge transfer associated with the HOMO-LUMO-excitation:
\be \label{charge}
    J[\phi_H,\phi_L] = \int_\Omega (x\cdot e) \Bigl( |\phi_L(x)|^2 - |\phi_H(x)|^2 \Bigr) dx
\ee
where $e$ is a given unit vector in $\R^3$ (the direction of charge transfer). 
Note that 
the integral $\int_\Omega (x\cdot e) |\chi|^2$ corresponds to the expected position of an electron with wavefunction $\chi$ along the $e$-axis, and so the integral in the definition of $J$ is the difference between the expected positions of LUMO and HOMO. Hence an excitation with maximal $J$ transfers electronic charge over a maximal distance. 
\\[2mm]
{\bf b) Overlap.} The spatial overlap between electron and hole can be measured by the functional
\be \label{over}
   J[\phi_H,\phi_L] = \int_\Omega |\phi_H|^2 |\phi_L|^2.  
\ee
Excitations with low $J$ have only a small amount of spatial overlap.
\\[2mm]
{\bf c) Lifetime.} The lifetime of excitations is of great interest in photovoltaics as a sufficient lifetime is needed for harvesting the photovoltaic current; estimating and controlling it in a simple manner is an important modelling challenge which we now address. 

The lifetime is governed, physically, by the subsequent time evolution of the system after the LUMO orbital has become occupied and the HOMO orbital has been vacated due to photon absorption. Thus to determine or control the lifetime one needs, ideally, to dynamically monitor (a model of) the subsequent time evolution of the excited system. 

As an -- analytically and computationally much more tractable -- substitute we propose the following ``time-infinitesimal'' way of quantitatively estimating and designing dynamic stability. Let $\phi$, $\phi_H$, $\phi_L$ denote, respectively, the KS orbitals \eqref{KS}, the HOMO \eqref{VHOMO}, and the LUMO \eqref{VLUMO}. Assume that the subsequent time evolution after excitation is governed by time-dependent density functional theory (TDDFT), with same map from density $\rho$ to exchange-correlation potential $v_{xc}(\rho)=\tfrac{d}{d\rho} e_{xc}(\rho)$ as the static theory (the latter simplification is known as the adiabatic local density approximation (ALDA)):
\begin{equation} \label{TDDFT}
    i\partial_t \phi'_i(\cdot,t) = h_{\rho'(\cdot, t)} \phi'_i(\cdot, t) \; (i=1,...,n),
    \;\;\; \rho'(\cdot, t) = 2\sum_{i=1}^n |\phi_i'(\cdot, t)|^2, \;\;\; \phi'_i(\cdot,0)=\phi'_i \; (i=1,...,n)
\end{equation}
with initial conditions given by the new orbitals after excitation, 
\begin{equation} \label{neworbs}
    \{ \phi'_1,...,\phi'_{n-1} \} \mbox{ any orthonormal basis of } 
    \{ \chi\in \mbox{Span}\,\{\phi_1,...,\phi_n\} \, : \, \langle \phi_H,\chi\rangle = 0\}, \;\;\; 
    \phi'_n = \phi_L.
\end{equation}
In terms of density matrices (mathematically: projectors onto the span of the occupied orbitals), the old occupied space span$\, \{\phi_1,...,\phi_n\}$ corresponds to the old density matrix $\gamma = \sum_{i=1}^n |\phi_i\rangle \langle\phi_i|$ while the new occupied space span$\, \{\phi'_1,...,\phi'_n\}$ corresponds to the new density matrix
\be \label{gamma'}
    \gamma' = \sum_{i=1}^n |\phi_i\rangle \langle\phi_i| - |\phi_H\rangle \langle \phi_H| + |\phi_L\rangle\langle\phi_L| = \gamma - |\phi_H\rangle \langle \phi_H| + |\phi_L\rangle\langle\phi_L|.
\ee
(As customary, we write $|\chi\rangle \langle \chi|$ for the orthogonal projector of $L^2(\Omega)$ onto the span of a normalized element $\chi$.) Note that $\gamma'$ depends neither on the choice of basis of old occupied space nor on the choice of basis of the orthogonal complement to $\phi_H$ in \eqref{neworbs}. In particular, the initial density $\rho'=\rho'(\cdot,0)$ in \eqref{TDDFT} is 
\be \label{rho'}
  \rho' = 2\sum_{i=1}^n|\phi_i|^2 - 2|\phi_H|^2 + 2|\phi_L|^2 = \rho - 2|\phi_H|^2 + 2|\phi_L|^2.
\ee
The failure of $\phi'$ to remain stationary up to unitary transformation, or equivalently to satisfy the stationary KS equations \eqref{KSeqs} with $h_\rho$ replaced by $h_{\rho'}$, is measured by the commutator $[h_{\rho'},\gamma']$. As a natural ``goal functional'' whose minimization promotes a long lifetime we thus propose the Hilbert-Schmidt norm of the commutator, 
\be \label{life}
   J[\phi,\phi_H,\phi_L] = \| [h_{\rho'},\gamma'] \|_{HS}^2 = \mbox{tr} \, \Bigl( [h_{\rho'},\gamma']^* [h_{\rho'},\gamma']\Bigr), \; \mbox{ with }\gamma' \mbox{ as in \eqref{gamma'}}.
\ee
Here $( \, )^*$ denotes the adjoint and $\mbox{tr}$ the trace. We call \eqref{life} the {\it lifetime functional}. 
Note that it is a functional purely of the occupied KS orbitals $\phi$, the HOMO $\phi_H$, and the LUMO $\phi_L$. (For an -- analytically and numerically convenient -- more explicit expression see Section \ref{S:optex}.) By construction it vanishes if and only if $\gamma'$ is invariant under the time evolution \eqref{TDDFT}, or equivalently $\phi'$ is time-invariant up to a time-dependent unitary transformation. 

An important remark is that, due to \eqref{KSeqs} and \eqref{HOMOLUMOeqs}, the commutator between $\gamma'$ and the old Hamiltonian vanishes, i.e. $[h_\rho,\gamma']=0$; consequently the commutator in \eqref{life} satisfies 
\be \label{commdiff}
   [h_{\rho'},\gamma'] = [h_{\rho'}-h_\rho, \gamma'],
\ee
that is to say it depends only on the {\it difference} of the KS Hamiltonians before and after excitation. Moreover this difference only comes from the Hartree- and the exchange-correlation term, 
\be \label{hdiff}
   h_{\rho'}-h_\rho = \frac{1}{|\cdot |} * 2 \bigl(|\phi_L|^2 - |\phi_H|^2\bigr) + \bigl( v_{xc}(\rho')-v_{xc}(\rho)\bigr). 
\ee
Physically, the first term in \eqref{hdiff} is the {\it (nonlocal) dipole field of the electron-hole pair}; the second term is the {\it (local) change in exchange-correlation potential caused by the electron-hole pair}. 

To summarize: Minimizing the goal functional \eqref{life} promotes excitations with a long lifetime. 
\\[2mm]
d) {\bf Bandgap.} A first basic question is: which range of HOMO-LUMO bandgaps $\eps_H-\eps_L$ can be engineered by suitably arranging the nuclear charge distribution $\mu$? This can be answered by maximizing respectively minimizing the bandgap, 
\be \label{gap}
   J[\phi,\phi_H,\phi_L,\mu] = \eps_L - \eps_H = \langle \phi_L, h_\rho \phi_L\rangle - \langle \phi_H, h_\rho \phi_H\rangle = \calE_{\mu,\rho}[\phi_L] - \calE_{\mu,\rho}[\phi_H].
\ee
Note that here, due to the occurrence of the $\mu$-dependent external potential $v_{ext}$ in the excitation functional $\calE_{\mu,\rho}$, the goal functional depends explicitly on $\mu$. We also note a formal similarity of the bandgap functional and the charge transfer functional \eqref{charge}: both are differences of expectation values for HOMO and LUMO, with the observable given, respectively, by a single-particle energy or Hamiltonian (3rd expression in \eqref{gap}) and a position operator. 

If one wants to tune the bandgap to a prescribed target value $\eps_*$, one needs to minimize a goal functional which reaches its minimum when $\eps_L-\eps_H=\eps_*$, the perhaps simplest one being 
\be \label{designgap}
    J[\phi,\phi_H,\phi_L,\mu] = \Bigl| \calE_{\mu,\rho}[\phi_L] - \calE_{\mu,\rho}[\phi_H] - \eps_*\Bigr|^2. 
\ee
\section{Existence of optimal excitations} \label{S:optex}
Our goal in this section is to prove that all the optimal control problems introduced in the previous section possess solutions.
\subsection{Analytic set-up} 
a) exchange-correlation energy: for the analysis below it suffices to assume 
\be \label{excassptns}
   e_{xc} \, : \, [0,\infty)\to\R \mbox{ continuously differentiable}, \; 
   |e_{xc}(\rho)|\le c_{xc}(1 \! +\! \rho^p ), \; |v_{xc}(\rho)| \le c_{xc}(1\! +\! \rho^{p-1})
\ee
for some exponent $p$ with $1\le p < \tfrac{5}{3}$ and some constant $c_{xc}$ independent of $\rho$. These assumptions are satisfied for all variants of the LDA which are used in practical calculations. In particular, they hold for the prototypical Dirac exchange energy \eqref{dirac}, with $p=\tfrac{4}{3}$. The condition $p<\tfrac{5}{3}$ guarantees that the KS energy functional is bounded from below, and is sharp for that purpose (see the proof of Lemma \ref{L:en} below).
\\[1mm]
b) orbitals: we assume that the minimization with respect to orbitals $(\phi_1,...,\phi_n)$ in the KS variational principle \eqref{KS} is over the admissible set
\be \label{An}
  \calA = \{ (\phi_1,...,\phi_n ) \in (H^1_0(\Omega))^n \, : \, \langle \phi_i, \, \phi_j\rangle = \delta_{ij} \},
\ee
where $H^1_0(\Omega)$ denotes the usual Sobolev space of square-integrable functions with square-integrable derivative which vanish on the boundary $\partial\Omega$, and $\langle \phi_i, \, \phi_j\rangle = \int_\Omega \overline{\phi_i} \, \phi_j$ denotes the $L^2$ inner product. Analogously, the optimization with respect to candidate HOMO and LUMO orbitals $\phi_H$ and $\phi_L$ is over the following sets which depend on the collection $\phi = (\phi_1,...,\phi_n)$ of minimizing KS orbitals from \eqref{KS}: 
$$
   \calA^H_\phi = \Bigl\{ \phi_H \in \mbox{ Span} \, \{\phi_1,...,\phi_n \} \, : \, \langle \phi_H,\, \phi_H\rangle = 1 \Bigr\}
$$
and 
$$
   \calA^L_\phi = \Bigl\{ \phi_L \in H^1_0(\Omega) \, : \, \langle \phi_i,\phi_L\rangle = 0 \mbox{ for }i=1,...,n, \; \langle \phi_L, \, \phi_L\rangle = 1 \Bigr\}.
$$
The governing variational problems \eqref{KS}, \eqref{VHOMO}, \eqref{VLUMO} for the Kohn-Sham, HOMO, and LUMO orbitals (the {\it state equation}) can now be written in the compact form
\begin{equation} \label{stateeq}
    \phi \in \underset{\calA}{\mbox{argmin \,}} \calE_\mu, \;\;\;
    \phi_H \in \underset{\calA_\phi^H}{\mbox{argmax \,}} \calE_{\mu,\phi}, \;\;\;
    \phi_L \in \underset{\calA_\phi^L}{\mbox{argmin \,}} \calE_{\mu,\phi}.
\end{equation}
c) nuclear charge distribution (the {\it control field}): we assume that the minimization or maximization with respect to the nuclear charge distribution in the control problems \eqref{goal}--\eqref{stateequation} is over nonnegative Radon measures of total mass $2n$ with support in $\Omega_{nuc}$, i.e. over the admissible set
\be \label{mu}
   \calA_{nuc} = \{ \mu \in \calM(\Omega_{nuc}) \, : \, \mu \ge 0, \; \int_{\Omega_{nuc}} d\mu = 2n \}, 
\ee 
where $\calM(\Omega_{nuc})$ (or $\calM$ for short) denotes the space of signed Radon measures on $\Omega_{nuc}$, i.e. the dual of the space $C(\Omega_{nuc})$ of continuous functions on $\Omega_{nuc}$, with norm $||\mu||_{\calM} = \sup \{ \int  f \, d\mu \, : \, f\in C(\Omega_{nuc}), \, \sup |f| \le 1\}$. 

One might wish to impose additional restrictions on $\mu$ such as: sum of delta functions (atomic nuclei) with integer nuclear charge; positions of the nuclei restricted to a crystal lattice; only a few types of atoms allowed, e.g.~``A'' and ``B''. Imposing all these would reduce optimization over \eqref{mu} to optimization over a set of doping profiles. Our theoretical analysis allows to easily incorporate such restrictions, see Theorem \ref{T:goal}. 
\subsection{Existence of excitations} \label{S:exci}
We begin by recalling standard results on existence of minimizers and analytic properties of the KS functional \eqref{KS} (see, e.g., \cite{CancesEtAl} for the more difficult case of an unbounded domain), and state a useful list of explicit estimates and continuity properties on function spaces which account for the dependence on the mass norm $||\mu||_{\cal M}$ of the nuclear charge distribution, as needed later (in section \ref{S:optoexci}) to show existence of optimal excitations. We then discuss the excitation functional \eqref{EHL} in the same spirit, and prove rigorously the existence of HOMO-LUMO excitations for general (measure-valued) nuclear charge distributions. Throughout, $||u||_p$ denotes the $L^p$ norm $(\int_\Omega |u|^p)^{1/p}$, $||u||_\infty$ stands for the $L^\infty$ norm $\sup_{x\in\Omega}|u(x)|$, and $H^1(\Omega)$ denotes the space of square-integrable functions on $\Omega$ with square-integrable derivative. Recall also that weak* convergence in $\calM(\Omega_{nuc})$ corresponds to convergence of continuous observables, i.e. a sequence $(\mu^{(\nu)})$ converges weak* to $\mu$ if $\int f\, d\mu^{(\nu)}$ converges to $\int f\, d\mu$ for all $f\in C(\Omega_{nuc})$. 
\begin{lemma} \label{L:en} (Lower bounds and continuity properties of the KS energy functional) The terms in the KS energy functional \eqref{KS} have the following properties: \\[0.5mm]
a) $T[\phi]\ge \tfrac12 T[\phi] + \tfrac{1}{4 c_s^2}||\rho||_3$, where $c_s$ is the Sobolev constant in the inequality $||u||_6 \le c_s ||\nabla u||_2$ on $\R^3$, and
$\phi\mapsto T[\phi]$ is continuous and weakly lower semicontinuous on $H^1(\Omega)^n$. 
\\[0.5mm]
b) $V_{ext}[\rho]\ge - ||\mu||_{\calM} \sup_{a\in\R^3}||\tfrac{1}{|\cdot -a|}||_2 ||\rho||_1^{1/4} ||\rho||_3^{3/4}$, and 
$(\phi,\mu) \mapsto V_{ext}[\rho]$ is strong $\times$ weak* continuous on $L^4(\Omega )^n\times \calM.$ \\[0.5mm]
c) $J_H\ge 0$, and $\phi\mapsto J_H[\rho]$ is continuous on $(L^{12/5}(\R^3))^n$. \\[0.5mm]
d) $E_{xc}[\rho] \ge - c_{xc}\bigl( \mbox{vol}(\Omega) + ||\rho||_1^{(3-p)/2} ||\rho||_3^{3(p-1)/2} \bigr)$, where $p\in[1,\tfrac{5}{3})$ is the exponent from \eqref{excassptns}, and 
$\phi\mapsto E_{xc}[\rho]$ is continuous on $L^{2p}(\Omega)^n$. 

In particular, $\calE_\mu[\phi]$ is continuous and weakly lower semicontinuous on $H^1(\Omega)^n$ as a functional of $\phi$, weak* continuous on $\calM$ as a functional of $\mu$, and strong $\times$ weak* continuous and weak $\times$ weak* lower semicontinuous on $H^1(\Omega)^n\times\calM$ as a functional of $(\phi,\mu)$.   
\end{lemma}
{\bf Proof} The estimate in a) is immediate from the well known inequality $T[\phi]\ge ||\nabla\sqrt{\rho}||_2^2$ and the Sobolev inequality applied with $u=\sqrt{\rho}$, and the continuity assertions are standard. The estimate in b) follows by first using the duality between $\calM(\Omega)$ and the space $C_b(\Omega)$ of bounded continuous functions on $\Omega$ and then the Cauchy-Schwarz inequality,
$$
   V_{ext}[\rho] = \int \Bigl( \frac{1}{|\cdot |}*\rho \Bigr)d\mu \ge -||\mu||_\calM ||\tfrac{1}{|\cdot|}*\rho||_\infty \ge - ||\mu||_{\calM} \sup_{a\in\R^3} ||\tfrac{1}{|\cdot -a|}||_2 ||\rho||_2, 
$$
and finally estimating the $2$-norm of $\rho$ by the H\"older interpolation inequality
$$
    ||\rho||_p \le ||\rho||_q^\theta ||\rho||_r^{1-\theta} \mbox{ with }q\le p \le r\mbox{ and }\tfrac{1}{p} = 
    \tfrac{\theta}{q} + \tfrac{1-\theta}{r},
$$
taking $p=2$, $q=1$, $r=3$. The continuity follows in a similar manner: the map $\rho\mapsto \tfrac{1}{|\cdot |}*\rho$ is continuous from $L^2(\Omega)$ to $C_b(\Omega)$ because $||\rho-\rho'||_\infty\le \sup_{a\in\R^3}||\tfrac{1}{|\cdot |}||_2 ||\rho-\rho'||_2$; and the map $(\phi_1,...,\phi_n)\mapsto \rho$ is continuous from $L^{4}(\Omega)^n$ to $L^2(\Omega)$. Hence the map $\phi\mapsto \tfrac{1}{|\cdot|}*\rho$ is continuous from $L^4(\Omega)^n$ to $C_b(\Omega)$. Since $V_{ext}$ is the duality pairing between $\tfrac{1}{|\cdot|}*\rho$ and $\mu$, the continuity assertion in b) follows. The lower bound on $J_H$ is trivial, and the asserted continuity is immediate from the well known continuity of $J_H$ as a functional of $\rho$ on $L^{6/5}(\R^3)$. Finally we deal with $E_{xc}$. By assumption \eqref{excassptns}, $E_{xc} \ge - c_{xc}\bigl(\mbox{vol}(\Omega) + ||\rho||_p^p\bigr)$. Applying the H\"older interpolation inequality with $q=1$, $r=3$ yields $\theta=(3-p)/(2p)$, $1-\theta=3(p-1)/(2p)$, and the asserted lower bound follows. As regards continuity, by continuity of the map $\phi\to\rho$ from $L^{2p}(\Omega)^n$ to $L^p(\Omega)$ it suffices to check that $\rho\mapsto e_{xc}(\rho)$ is continuous from $L^p(\Omega)$ to $L^1(\Omega)$; but this is a standard consequence of the pointwise continuity of $e_{xc}$ as a function of $\rho$ and the growth bound in \eqref{excassptns}. 
\\[2mm]
We now turn to the excitation functional.
\begin{lemma}\label{L:exen} (Lower bounds and continuity properties of the excitation functional) The terms in the excitation functional \eqref{EHL} have the following properties: denoting $\rho_\chi=|\chi|^2$, \\[1mm]
a) $T[\chi]\ge \tfrac{1}{2}T[\chi] + \tfrac{1}{4c_s^2}||\rho_\chi||_3$, and
$$
   \chi\mapsto T \mbox{ is continuous and weakly lower semicontinuous on }H^1(\Omega).
$$
b) $\int v_{ext} \rho_\chi \ge - ||\mu||_\calM \sup_{a\in\R^3}||\tfrac{1}{|\cdot - a|} ||_2 ||\rho_\chi||_1^{1/4} ||\rho_\chi||_3^{3/4}$, and 
$$
    (\chi,\mu)\mapsto \int v_{ext} \rho_\chi \mbox{ is strong $\times $ weak* continuous on }L^4(\Omega) \times \calM.
$$
c) $\int \bigl(\tfrac{1}{|\cdot |}*\rho_\chi)\rho \ge 0$, and
$$
    (\phi,\chi)\mapsto \int \bigl(\tfrac{1}{|\cdot |}*\rho_\chi\bigr) \rho \mbox{ is continuous on }
    L^{12/5}(\Omega)^{n+1}.
$$
d) $\int v_{xc}(\rho)\, \rho_\chi \ge - c_{xc}\bigl(||\rho_\chi||_1 + ||\rho||_p^{p-1} ||\rho_\chi||_1^{(3-p)/(2p)} ||\rho_\chi||_3^{3(p-1)/(2p)} \bigr)$, and 
$$
    (\phi,\chi) \mapsto \int v_{xc}(\rho)\, \rho_\chi \mbox{ is continuous on }
    L^{2p}(\Omega)^{n+1}.
$$
In particular, $(\phi,\chi,\mu)\mapsto \calE_{\mu,\phi}[\chi]$ is weak $\times$ strong $\times$ weak* continuous and weak $\times$ weak $\times$ weak* lower semicontinuous on $(H^1(\Omega))^n \times H^1(\Omega) \times \calM$.
\end{lemma}
{\bf Proof} The results for $T$, $V_{ext}$, and $\int (\tfrac{1}{|\cdot|}*\rho_\chi)\rho$ follow analogously to those for $T$, $V_{ext}$ and $J_H$ in Lemma \ref{L:en}. To deal with the exchange-correlation term is a little more work. By the bound on $v_{xc}$ in \eqref{excassptns} and H\"older's inequality we have 
$$
   \int v_{xc}(\rho) \, \rho_\chi \ge - c_{xc} \int\bigl(1 + \rho^{p-1}\bigr) \rho_\chi \ge - c_{xc}\bigl( ||\rho_\chi||_1 + \underbrace{||\rho^{p-1}||_{p'}}_{=||\rho||_p^{p-1}} ||\rho_\chi||_p\bigr) \mbox{ with }p'=\tfrac{p}{p-1}.
$$ 
The asserted bound now follows from the H\"older interpolation inequality with $q=1$, $r=3$. The continuity is clear when $p=1$, so let $p>1$. By the pointwise continuity of $\rho\mapsto v_{xc}(\rho)$ and the bound on $v_{xc}$ in \eqref{excassptns}, $\rho\mapsto v_{xc}(\rho)$ is continuous from $L^{q(p-1)}$ to $L^q$ whenever $\min\{1,\tfrac{1}{p-1}\}\le q<\infty$. Taking $q=p'=p/(p-1)$, we obtain continuity from $L^p$ to $L^{p'}$. Hence, since $L^p$ and $L^{p'}$ are in duality, the map $(\rho,\rho_\chi)\mapsto \int v_{xc}(\rho) \, \rho_\chi$ is continuous on $L^p(\Omega)\times L^p(\Omega)$, yielding the asserted continuity. 
\\[2mm]
Next we collect basic analytic properties of the admissible sets of trial functions in the variational problems \eqref{stateeq}. Thanks to the compact embedding of $H^1(\Omega)$ into $L^2(\Omega)$ and the obvious $L^2$-continuity of the constraints appearing in the definitions of the sets we have:
\begin{lemma} \label{L:constraint} (Analytic properties of the constraints) \\
a) $\calA$ is weakly closed in $H^1(\Omega)^n$. \\
b) For any $\phi\in\calA$, the sets $\calA_\phi^H$ and $\calA_\phi^L$ are weakly closed in $H^1(\Omega)$. 
\end{lemma}
It is now a straightforward matter to recover the following well known fact:
\begin{lemma} \label{L:existence} (Existence of KS ground states) For any nuclear charge distribution $\mu\in\calA_{nuc}$, there exists a minimizer $\phi=(\phi_1,...,\phi_n)$ of the KS energy functional $\calE_{\mu}$ on the admissble set $\calA$.
\end{lemma}
Moreover we straightforwardly infer the following new result:
\begin{theorem} \label{T:existence} (Existence of HOMO-LUMO excitations) For any nuclear charge distribution $\mu\in\calA_{nuc}$, and any set of orbitals $\phi=(\phi_1,...,\phi_n)\in\calA$, the excitation functional \eqref{EHL} possesses a maximizer $\phi_H$ on $\calA_\phi^H$ (i.e., a HOMO)  and a minimizer $\phi_L$ on $\calA_\phi^L$ (i.e., a LUMO). 
\end{theorem}
{\bf Proof of Lemma \ref{L:existence}} By the estimates in Lemma \ref{L:en} together with the fact that the exponent $\tfrac{3(p-1)}{2}$ of $||\rho||_3$ in d) is $<1$ (under the sharp assumption on the exchange-correlation energy made in \eqref{excassptns} that $p<\tfrac{5}{3}$) and that $||\mu||_\calM=\int d\mu = 2n$ (see \eqref{mu}), $\calE_\mu$ is bounded from below and there exists a constant $C(n,p,c_{xc},\Omega)$ such that, whenever $(\phi,\mu)\in\calA\times\calA_{nuc}$ and $\calE_\mu[\phi]\le \inf_\calA \calE_\mu + 1$, then 
\be \label{upbd}
   ||\phi ||_{H^1} \le C(n,p,c_{xc},\Omega ).
\ee 
The assertion now follows from a standard lower semicontinuity / compactness argument: any minimizing sequence possesses a subsequence converging weakly in $H^1(\Omega)^n$; 
and the weak limit belongs, by Lemma \ref{L:constraint}, to the admissible set, and is a minimizer, 
by the weak lower semicontinuity result at the end of Lemma \ref{L:en}.
\\[2mm]
{\bf Proof of Theorem \ref{T:existence}} Existence of a HOMO is straightforward since $\calE_{\mu,\phi}$ is continuous on $H^1(\Omega)$ (see Lemma \ref{L:exen}) and the admissible set $\calA_{\phi}^H$ -- being a closed bounded subset of a finite-dimensional subspace -- is compact. Existence of a LUMO follows by arguing similarly to the proof of Lemma \ref{L:existence}. For completeness and future reference we include the details. By the bounds in Lemma \ref{L:exen} and the fact that the exponent $\tfrac{3(p-1)}{2p}$ is $<1$ thanks to the sharp assumption $p<\tfrac{5}{3}$ in \eqref{excassptns}, the excitation functional \eqref{EHL} is bounded from below. Moreover there exists a constant $D(n,p,c_{xc},\Omega,C)$ depending only on $n$, $p$, $c_{xc}$, $\Omega$, and an upper bound $C$ on the $H^1$ norm of $\phi$ such that whenever $\mu\in\calA_{nuc}$, $\phi\in\calA$ with $||\phi||_{H^1}\le C$, and $\calE_{\mu,\phi}[\chi]\le \inf_{\calA_{\phi}^L}\calE_{\mu,\phi}+1$, then 
\be \label{upbdchi}
   ||\chi||_{H^1} \le D(n,p,c_{xc},\Omega, C).
\ee
In particular any minimizing sequence of $\calE_{\mu,\phi}$ on $\calA_{\phi}^L$ possesses a subsequence which is weakly convergent in $H^1(\Omega)$. The lower semicontinuity result for $\calE_{\mu,\phi}$ in Lemma \ref{L:exen} together with the closedness result on $\calA_\phi^L$ in Lemma \ref{L:constraint} implies that the weak limit is a minimizer. 
\subsection{The set of HOMO-LUMO excitations} 
The main technical idea underlying existence of optimal excitations lies in introducing and analyzing the set of {\it all} $\,$HOMO-LUMO excitations which can be achieved by some (measure-valued) distribution of nuclear charge. This idea allows us to overcome the difficulty that the goal functionals, when viewed as functionals only of the nuclear charge distribution $\mu$ via the map $\mu\mapsto (\phi,\phi_H,\phi_L)$ defined by \eqref{stateeq}, do not appear to have any useful continuity properties due to the lack of convexity and uniqueness in the problems \eqref{stateeq}.
\begin{lemma} \label{L:state} (Analytic properties of the set of HOMO-LUMO excitations) \\
The joint solution set to the governing variational principles \eqref{stateeq} for occupied KS orbitals, HOMO, and LUMO parametrized by the set of nuclear charge distributions $\mu$,
$$
   \calB = \{ (\phi,\phi_H,\phi_L,\mu) \, : \, \mu\in\calA_{nuc}, \, \eqref{stateeq} \},
$$
has the following properties: \\[1mm]
a) It is weak $\times$ weak $\times$ weak $\times$ weak* closed in 
$H^1(\Omega)^n\times H^1(\Omega)\times H^1(\Omega)\times\calM$. \\[1mm]
b) It is strong $\times$ strong $\times$ strong $\times$ weak* compact in $H^1(\Omega)^n\times H^1(\Omega)\times H^1(\Omega)\times\calM$.
\end{lemma}
{\bf Proof} This statement is more subtle than the results in the previous section. We begin with a). 
Let $(\phinu,\phinuH, \phinuL)\weakly(\phi,\phi_H,\phi_L)$ in  $H^1(\Omega)^{n+2}$ and $\munu\weakly^*\mu$ in $\calM$, where here and below the halfarrows $\weakly$ and $\weakly^*$ stand for weak and weak* convergence. We need to check that 
(i) $\mu\in\calA_{nuc}$, 
(ii) $\phi\in\mbox{argmin}_\calA \calE_\mu$, 
(iii) $\phi_H\in\mbox{argmax}_{\calA_\phi^H}\calE_{\mu,\phi}$, 
(iv) $\phi_L\in\mbox{argmin}_{\calA_\phi^L}\calE_{\mu,\phi}$.

(i) is clear from the fact that the $\munu$ are supported in the compact set $\Omega_{nuc}$. 

(ii) does not follow from the weak lower semicontinuity of $\calE_\mu$ on $H^1(\Omega)^n$, since the nuclear charge distribution is now also varying. Using first the variational property of $\phinu$ and then the weak* continuity of the map $\mu\mapsto \calE_\mu[\psi]$ for fixed $\psi\in\calA$ (see Lemma \ref{L:en}), we have
\be \label{firstineq}
    \calE_{\munu}[\phinu] \le \calE_{\munu}[\psi] \to \calE_\mu[\psi]
\ee
and hence $\limsup_{\nu\to\infty} \calE_{\munu}[\phinu]\le \calE_\mu[\psi]$. Since $\psi\in\calA$ was arbitrary, it follows that
\be \label{upbdii}
       \limsup_{\nu\to\infty} \calE_{\munu}[\phinu ]\le \inf_\calA \calE_\mu.
\ee
On the other hand, by the weak $\times$ weak* lower semicontinuity of $(\phi,\mu)\mapsto \calE_\mu[\phi]$ on $H^1(\Omega)^n\times\calM$ we have 
\be \label{lowbdii}
       \liminf_{\nu\to\infty} \calE_{\munu}[\phinu ]\ge \calE_\mu[\phi].
\ee  
Inequalities \eqref{upbdii}, \eqref{lowbdii} together with the fact that by Lemma \ref{L:constraint} $\phi\in\calA$ (i.e., $\phi$ is an admissible trial function in the variational problem in \eqref{upbdii}) yields (ii). 

Moreover we conclude that $\calE_{\munu}[\phinu]\to \calE_\mu[\phi]$. But the functional $\calE_\mu$ consists of four contributions, of which $V_{ext}$, $J_H$, and $E_{xc}$ have the property that their values for the sequence $(\phinu,\munu)$ converge to those for the limit $(\phi,\mu)$ (see Lemma \ref{L:en}). Hence, importantly, the remaining term $T$ must satisfy 
$T[\phinu] \to T[\phi]$.    
Consequently $\nabla\phinu\weakly\nabla\phi$ in $L^2$ and $||\nabla\phinu||_2\to||\nabla\phi||_2$. These two statements together imply $\nabla\phinu\to\nabla\phi$ strongly in $L^2$, and hence
$\phinu\to\phi$ in $H^1(\Omega)^n$. This fact will be useful later. 

To deal with (iii) and (iv) is more difficult, since the starting point of the above argument -- the first inequality in \eqref{firstineq} -- now {\it fails} as the HOMO and LUMO orbitals $\phinuH$ and $\phinuL$ do not belong to universal but rather to $\phinu$-dependent sets, and hence candidate orbitals $\chi$ in the variational principles for the limiting HOMO and LUMO orbitals $\phi_H$ and $\phi_L$ are not admissible trial functions in the variational principle for the approximating orbitals. The idea, then, is to bring into play the $L^2$ projector 
$
          \gamma_{\phinu}\chi := \sum_{i=1}^n \langle \phinui,\chi\rangle \phinui.
$ 
Clearly, for fixed $\chi\in L^2(\Omega)$ the map $\psi\mapsto \gamma_\psi\chi$ is strongly continuous from $H^1(\Omega)^n$ to $H^1(\Omega)$; in particular, 
\be \label{star}
    \gamma_{\phinu} \chi \to \gamma_\phi \chi \mbox{ in }H^1(\Omega), \;\;\; 
    ||\gamma_\phinu \chi ||_2 \to ||\gamma_\phi \chi ||_2. 
\ee
Now let $\chi\in\calA_\phi^H$. It follows that $\gamma_\phi\chi = \chi$ and $||\chi||_2=1$. Hence by \eqref{star} $||\gamma_\phinu \chi||_2>0$ for all sufficiently large $\nu$, and -- by the variational principle for the HOMO $\phinuH$ -- 
\be \label{firstineqiii}
   \calE_{\munu,\phinu}[\phinuH ] \ge \calE_{\munu,\phinu} \Bigl[ \frac{\gamma_\phinu  \chi}{||\gamma_\phinu \chi ||_2} \Bigr] = \frac{1}{||\gamma_\phinu \chi||_2^2} \calE_{\munu,\phinu} [\gamma_\phinu \chi ] \; \to \; 1 \cdot \calE_{\mu,\phi}[\chi ],
\ee 
with the above convergence being due to the continuity property of the map $(\phi,\chi,\mu)\mapsto \calE_{\mu,\phi}[\chi]$ proved in Lemma \ref{L:exen}. Since \eqref{firstineqiii} is true for all $\chi\in\calA_\phi^H$, it follows that
\be \label{lowbdiii}
    \liminf_{\nu\to\infty} \calE_{\munu,\phinu} [\phinuH] \ge \sup_{\chi\in\calA_\phi^H} \calE_{\mu,\phi}[\chi ], 
\ee
the right hand side being the limiting HOMO eigenvalue $\eps_H$. 
On the other hand, since $\phinuH\in\calA_\phinu^H$ we have
$
   \phinuH = \gamma_{\phinu} \phinuH = \sum_{i=1}^n \langle \phinui, \phinuH\rangle \phinui,
$
and so, by the weak convergence of $\phinuH$ in $H^1$ and the strong convergence of $\phinu$ in $H^1$, $\phinuH \to \sum_{i=1}^n \langle \phi_i,\phi_H\rangle \phi_i = \gamma_\phi \phi_H$ strongly in $H^1(\Omega )$. 
Hence $\phi_H=\gamma_\phi\phi_H$, that is to say $\phi_H\in\calA_\phi$. Moreover by the weak $\times$ strong $\times$ weak* continuity of $(\phi,\chi,\mu)\mapsto \calE_{\mu,\phi}[\chi]$ we have 
\be \label{upbdiii}
   \calE_{\munu,\phinu}[\phinuH ] \to \calE_{\mu,\phi} [\phi_H ].
\ee
The upper and lower bounds \eqref{lowbdiii}, \eqref{upbdiii} together imply (iii). 

(iv): Let $\chi\in\calA_\phi^L$. Thus $\gamma_\phi \chi = 0$, $||\chi||_2=1$. Hence by \eqref{star}, 
$(I-\gamma_\phinu)\chi \to \chi$ in $H^1(\Omega)$, 
$||(I-\gamma_\phinu) \chi||_2\to 1$. By the variational principle for the LUMO $\phinuH$, 
$$
     \calE_{\munu,\phinu}[\phinuL] \le \calE_{\munu,\phinu} \Bigl[ \frac{(I-\gamma_\phinu)\chi}{||(I-\gamma_\phinu)\chi||_2} \Bigr] = 
     \frac{1}{||(I-\gamma_\phinu)\chi||_2^2} \calE_{\munu,\phinu}[(I-\gamma_\phinu)\chi] \to 1 \cdot \calE_{\mu,\phi}[\chi],
$$
the convergence above following as in \eqref{firstineqiii}. Minimization over $\chi\in\calA_\phi^L$ yields 
\be \label{upbdiv}
   \limsup_{\nu\to\infty} \calE_{\munu,\phinu}[\phinuL ] \le \inf_{\chi\in\calA_\phi^L} \calE_{\mu,\phi}[\chi ],
\ee 
the right hand side being the limiting LUMO eigenvalue $\eps_L$. The proof of the complementing lower bound is {\it different} as for $\phi_H$, since -- at this point -- the $\phinuL$ are {\it not} known to converge strongly in $H^1$ but only weakly. On the other hand we now only need to bound the excitation energy of the limit from above, for which we can use the weak $\times$ weak $\times$ weak* lower semicontinuity of $(\phi,\chi,\mu)\mapsto\calE_{\mu,\phi}[\chi]$:
\be \label{lowbdiv}
    \calE_{\mu,\phi}[\phi_L ] \le \liminf_{\nu\to\infty} \calE_{\munu,\phinu }[\phinuL ].
\ee
Together with \eqref{upbdiv} this establishes (iv), completing the proof of a). Moreover we conclude that $\calE_{\munu,\phinu}[\phinuL] \to \calE_{\mu,\phi}[\phi_L]$. This implies that $\phinuL\to\phi_L$ strongly in $H^1(\Omega)$, by arguing analogously to the proof of the strong convergence of $\phinu$ in $H^1$ but using Lemma \ref{L:exen} instead of Lemma \ref{L:en}. In total we have shown that $(\phinu,\phinuH,\phinuL)$ converges strongly in $H^1(\Omega)^{n+2}$.

Thus to establish b), by the Banach-Alaoglu theorem it suffices to show that any sequence $(\phinu,\phinuH, \phinuL, \munu)\in\calB$ is bounded in $H^1(\Omega)^n\times H^1(\Omega)\times H^1(\Omega)\times\calM$. For $\phinu$, $\phinuL$ and $\munu$ this follows from, respectively, \eqref{upbd}, \eqref{upbdchi}, and $||\munu||_\calM = \int\munu = 2n$. Finally, boundedness of $\phinuH$ follows straightforwardly from the boundedness of $\phinu$.
%
%
%
\subsection{Optimal excitations} \label{S:optoexci}
We are finally in a position to show well-posedness of all the novel optimal control problems introduced in Section \ref{S:OCPs}. In the sequel, the goal functionals are always considered as functionals of $\phi$ (occupied KS orbitals), $\phi_H$ (HOMO), $\phi_L$ (LUMO), $\mu$ (nuclear charge distribution), even though some of them depend only on HOMO and LUMO. 
\begin{theorem} \label{T:goal} (Existence of optimal excitations) The optimal control problem to maximize or minimize \eqref{goal} over $(\phi,\phi_H,\phi_L,\mu)\in H^1(\Omega)^{n+2}\times\calA_{nuc}$ subject to the constraint \eqref{stateeq} possesses a solution when the goal functional is any of the functionals \eqref{charge}, \eqref{over}, \eqref{life}, \eqref{gap}, \eqref{designgap}. More generally, it possesses a solution whenever the goal functional is continuous on $\calB$ with respect to strong $\times$ weak* convergence in $H^1(\Omega)^{n+2}\times\calM$ and $\calA_{nuc}$ is replaced by any weak* closed subset.  
\end{theorem}
{\bf Proof} The main work has already been carried out in Lemma \ref{L:state} b), which immediately implies the second part of the theorem, and reduces the first part to checking the required continuity of the goal functionals, i.e. to proving:
\begin{lemma} \label{L:goal} The charge transfer functional \eqref{charge}, the overlap functional \eqref{over}, the lifetime functional \eqref{life}, and the bandgap functionals \eqref{gap} and \eqref{designgap} are continuous on $\calB$ with respect to strong $\times$ weak* convergence in $H^1(\Omega)^{n+2}\times \calM$. 
\end{lemma}
Importantly, Lemma \ref{L:state} b) provides strong rather than just weak compactness in $H^1$ for the HOMO and LUMO orbitals; we note that weak compactness would not be enough as, e.g., the bandgap does not appear to be continuous or upper semicontinuous or lower semicontinuous with respect to weak $H^1$ convergence.

{\bf Proof of Lemma \ref{L:goal}} The functionals \eqref{charge} and \eqref{over} only depend on $\phi_H$ and $\phi_L$, and are obviously continuous on, respectively, $L^2\times L^2$ and $L^4\times L^4$. The bandgap \eqref{gap} is, by Lemma \ref{L:exen}, weak $\times$ strong $\times$ strong $\times$ weak* continuous on $(H^1)^{n}\times H^1 \times H^1 \times \calM$, establishing in particular the asserted continuity of \eqref{gap} and \eqref{designgap}. 

It remains to look at the lifetime functional \eqref{life}. A long but rather elementary calculation yields the following more explicit expression which only involves a sum over {\it occupied} orbitals in the excited state:  
\be \label{life'}
    \mbox{tr}\, \Bigl( [h',\gamma']^*[h',\gamma']\Bigr) = 2 \sum_{\phi \in \{ \phi_i-\langle\phi_H,\phi_i\rangle\phi_H\}_{i=1}^n \cup \{ \phi_L \} } \| (I-\gamma') (h_{\rho'}-h_\rho) \phi \|_2^2, \;\;\; \gamma' \mbox{ as in }\eqref{gamma'}.   
\ee
Note that the operator $I-\gamma'$ appearing on the right hand side is the projector onto {\it unoccupied} space. 
(The above expression has the additional virtue of being invariant under the choice of basis $\phi_1,..,\phi_n$ of the occupied orbital space in the ground state. If the HOMO $\phi_H$ is the $n^{th}$ KS orbital $\phi_n$, then the corresponding state $\phi_n-\langle\phi_H,\phi_n\rangle \phi_H$ is zero, as is its contribution to the above sum, and so the sum reduces to a sum over $\phi\in\{\phi_1,...,\phi_{n-1},\phi_L\}$.) To derive \eqref{life'}, one starts by using \eqref{commdiff} and expands the trace in the form $\mbox{tr}(A^* A) = \sum_{\ell,k=1}^\infty |\langle \phi_\ell, A\phi_k\rangle|^2$ in an ONB $\{\phi_1,...,\phi_n, \phi_L, \phi_{n+2}, \phi_{n+3}, ...\}$, the first $(n+1)$ basis functions of which consist of the occupied KS orbitals in the ground state and the LUMO. One then uses the explicit form of $\gamma'$ to evaluate the individual terms.  

We now analyze continuity of the right hand side of \eqref{life'}, using expression \eqref{hdiff} for the difference of the KS Hamiltonians. In the proof of Lemma \ref{L:en} we showed that $(\phi_H,\phi_L)\mapsto \tfrac{1}{|\cdot |}*2(|\phi_H|^2-|\phi_L|^2)=:\tilde{v}$ is continuous from $L^4\times L^4$ to $C_b$, and so $(\phi_H,\phi_L,\chi)\mapsto \tilde{v}\,\chi$ is continuous from $L^4\times L^4 \times L^2$ to $L^2$. Moreover, since by \eqref{excassptns} and the fact that $p<5/3$ we have $|v_{xc}(\rho)|\le c_{xc}(1+\rho^{2/3})$, the map $(\rho,\rho_{\phi_H},\rho_{\phi_L})\mapsto v_{xc}(\rho')- v_{xc}(\rho)$ is continuous from $(L^{7/3})^3$ to $L^{7/2}$, whence $(\phi,\phi_H,\phi_L,\chi)\mapsto (v_{xc}(\rho')-v_{xc}(\rho))\chi$ is continuous from $(L^{14/3})^n\times (L^{14/3})^3$ to $L^2$ (note that if $w\in L^{7/2}$ and $\chi\in L^{14/3}$, then the product $w\,\chi$ is in $L^2$, since $\tfrac{1}{7/2}+\tfrac{1}{14/3}=\tfrac{1}{2}$). Finally, consider the projector $(I-\gamma')$ onto unoccupied space. It is clear from the expression \eqref{gamma'} that $(\phi,\phi_H,\phi_L,\chi)\mapsto (I-\gamma')\chi$ is continuous from $(L^2)^{n+3}$ to $L^2$. In total, $(\phi,\phi_H,\phi_L,\chi) \mapsto (I-\gamma') (h_{\rho'}-h_\rho) \chi$ is continuous from $L^{14/3}(\Omega )^{n+3}$ to $L^2(\Omega )$.
Finally we need to look at the orbitals $\chi$ over which the sum in \eqref{life'} runs. The maps $(\phi,\phi_H,\phi_L)\mapsto \phi_i - \langle\phi_H,\phi_i\rangle\phi_H$ are obviously continuous from $(L^{14/3})^{n+2}$ to $L^{14/3}$. Thus \eqref{life'} as a functional of $(\phi,\phi_H,\phi_L)$ is continuous on $(L^{14/3})^{n+2}$, completing the proof of Lemma \ref{L:goal} and Theorem \ref{T:goal}.                                                                                             
\section{Numerical results}
\label{sec:numerics}
In order to illustrate the excitation properties achievable via optimal control, we investigate the following 1D model problem. We consider a 1D nanocrystal consisting of 20 atoms with fixed equidistant positions and 120 electrons, but allow the nuclear charges of the atoms to change while preserving a total nuclear charge of 120 (this value is chosen so that the undoped chain corresponds to pure carbon). 

More precisely, the nuclear positions are fixed at $R_1,...,R_{20}$ = $-9.5,\, -8.5, ..., \, 8.5, \, 9.5$ and the admissible nuclear charge distributions are taken to be sums of sharply peaked Gaussians of width $\sigma$, 
\be \label{mu1D}
   \mu(x) = \sum_{\alpha=1}^{20} Z_\alpha \frac{\exp(-\frac{(x-R_\alpha)^2}{2\sigma^2})}{\sqrt{2\pi\sigma^2}},
\ee
with integer nuclear charges $Z_1,...,Z_{20}\in \{3,...,9\}$ satisfying $\sum_{\alpha=1}^{20}Z_\alpha=120$. These bounds correspond to Li to F, and exclude noble gas atoms.
We note that there are $\sim 7^{20}$ configurations (ignoring the - in terms of order of magnitude negligible - constraint on the total charge). 

As regards the modelling of the electronic structure within Kohn-Sham DFT, to simplify computations and because of the lack of simple exchange-correlation functionals in 1D we drop exchange-correlation contributions, and employ the 1D Kohn-Sham Hamiltonian (see \eqref{KSeqs} for notation)
\be \label{ham1D}
   h_\phi = -\frac12 \, \frac{d^2}{dx^2} + v_{ext} + v_H,
\ee
in the spatial domain $[-10,10]$. Moreover, as the Coulomb potential $v(x)=1/|x|$ is not integrable in 1D, we replace it by the effective longitudinal Coulomb potential in a thin wire of diameter $d$
\cite{BSCA03, CF15},
\be \label{coul1D}
   v_d(x) = \frac{\sqrt{\pi}}{2d} \exp \left( \frac{x^2}{4d^2} \right) \mathrm{erfc} \left( \frac{x}{2d} \right).
\ee
As described in Section \ref{S:KS} we restrict ourselves to closed-shell electron configurations, so in the ground state the lowest $60$ orbitals $\phi_i$ are doubly occupied, HOMO equals orbital $60$, and LUMO orbital $61$. The adaptation of the goal functionals to the above 1D setting is straightforward; for instance, the charge transfer functional \eqref{charge} becomes
\be \label{charge1D}
   J[\phi_H,\phi_L] = \int_{-10}^{10} x\Bigl( |\phi_{61}(x)|^2 - |\phi_{60}(x)|^2 \Bigr) dx. 
\ee

In our simulations we chose the parameters in \eqref{mu1D} and \eqref{charge1D} to be $d=0.01$ and $\sigma^2=1/2000$, and for the discretization of space and time we employ grids of meshsize $\Delta x=0.01$ and $\Delta t = 0.002$. 
\begin{figure}
	\centering
	\hspace*{-.1525\textwidth} 
    \includegraphics[width=0.37\textwidth]{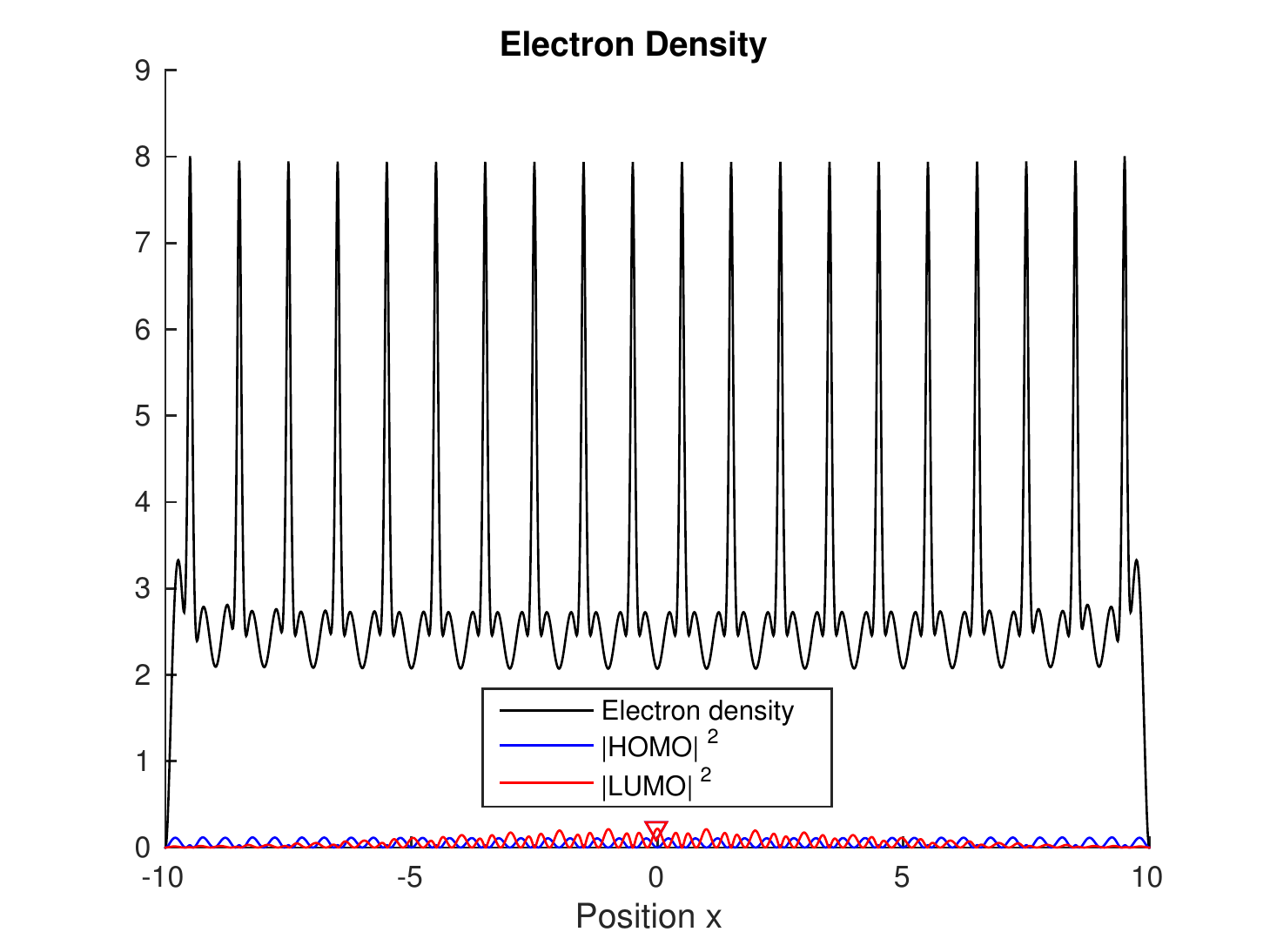}
	\hspace*{-.04\textwidth} \includegraphics[width=0.37\textwidth]{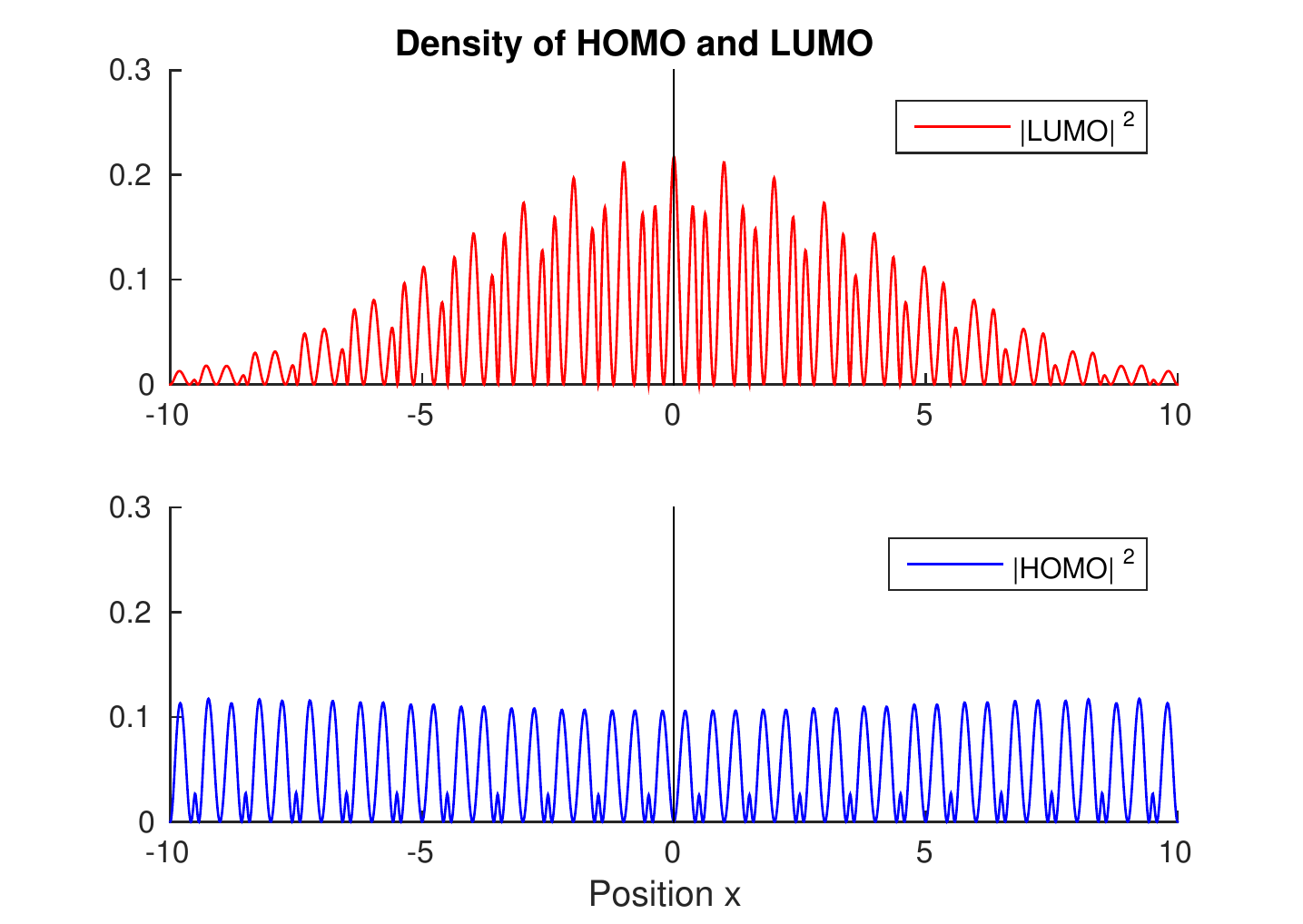}
	\hspace*{-.04\textwidth} \includegraphics[width=0.37\textwidth]{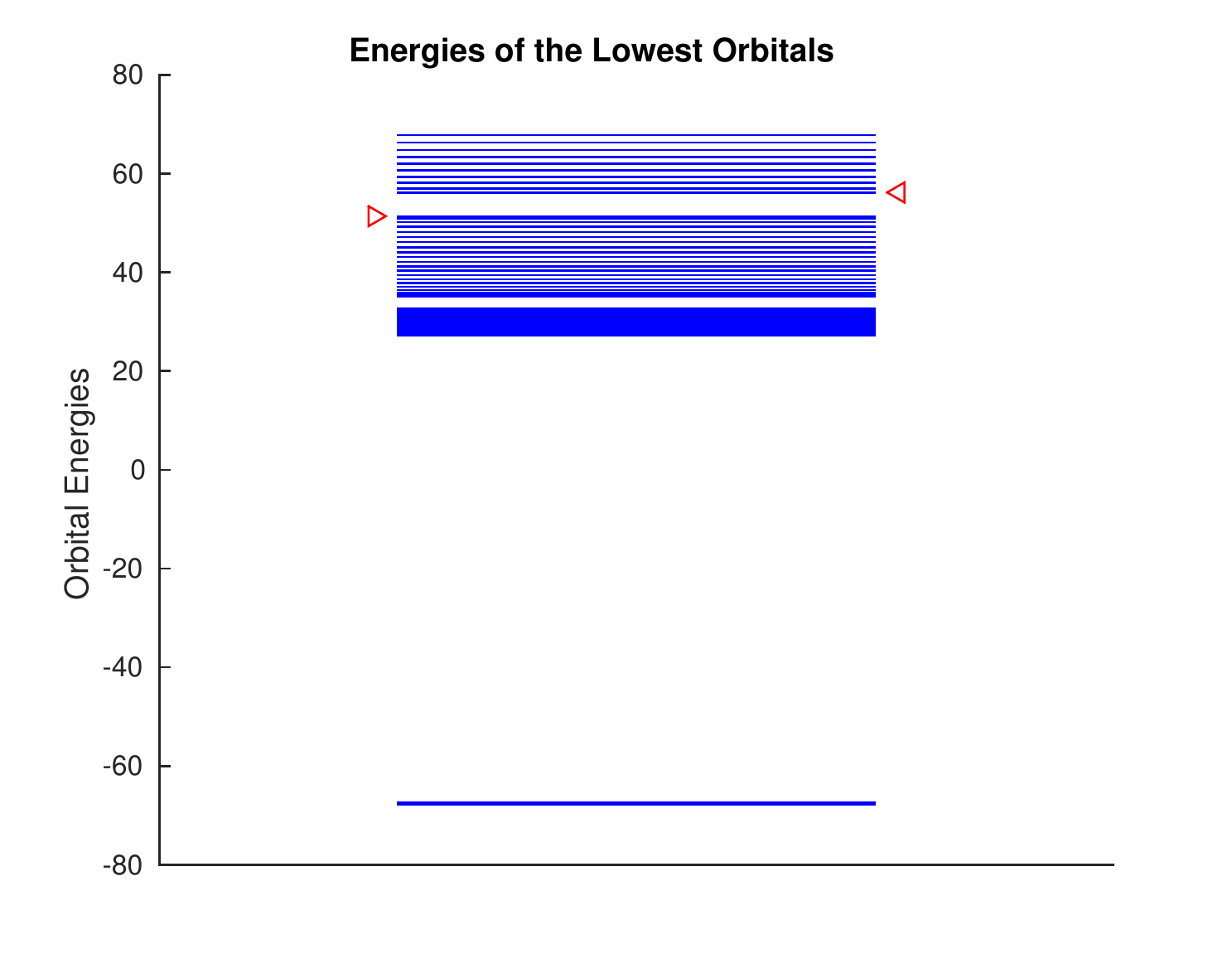} \hspace*{-.15\textwidth}
	\caption{Pure carbon chain. Left and center: The electron density of the ground state and the density of HOMO and LUMO. The center of mass of both HOMO and LUMO (vertical lines) is located at the midpoint of the chain, so the excitation produces no net charge transfer. 
Right: The single-particle energy levels (eigenvalues of the KS Hamiltonian (4.2) with ground state orbitals $\phi$), showing a typical band structure and a pronounced bandgap between HOMO and LUMO (small triangles).}
	\label{figcarbon}
\end{figure}
\vspace*{2mm}

\noindent
{\bf Optimization algorithm.}
Our approach to find atomic configurations which (approximately) optimize our various goal  functionals is a genetic algorithm based on randomly chosen directions $h \in \{-1,0,1\}^{20}$ in configuration space which we use to update the vector $Z=(Z_1,...,Z_{20})$ of nuclear charges to $Z + h$. Similarly to stochastic optimization techniques like simulated annealing, we carry out ``large'' steps $h$ at the beginning and gradually reduce the (Euclidean) length of $h$ to ``small'' steps, while at the same time increasing the number of search directions in order to allow for a more detailed exploration of the configuration space close to a possible local extremum. 
More precisely, our algorithm is as follows:
\begin{itemize}
\item Start from the pure carbon chain, $Z(0)=(6,6,...,6)$. 
\item In the i$^{th}$ optimization step, start from the previous optimum $Z(i-1)$; generate $n(i)$ random increments $h\in\{-1,0,1\}^{20}$ with probability $P(h_j\! =\!\pm 1)=p(i)$, zero total charge, and $Z+h\in\{3,...,9\}^{20}$; pick the best increment, i.e. $Z(i)=Z(i-1)+t_*h_*$ with $t_*\in\{-1,0,1\}$ and $h_*$ giving the best value for $J[Z+t_*h_*]$. 
\item Use geometrically decreasing probabilities $p(i)=p(1)\cdot 2^{-i}$ for components of the increment to be nonzero, and a geometrically increasing number $n(i)=n(1)\cdot 2^i$ of search directions. 
\end{itemize}
In our numerical results we took $p(1)=\tfrac{1}{3}$, $n(1)=10$, and $4$ iteration steps. Note that in the $4^{th}$ step, the probability of any entry $h_j$ being nonzero is only $\tfrac{1}{12}$; but since at least two nonzero entries are needed because of charge conservation, this already corresponds to maximal sparsity -- or equivalently minimal size -- of increment vectors.   

We stress that our optimization algorithm is stochastic. Therefore one obtains different final configurations (or doping profiles) $Z_*$ in each run of the algorithm. However the associated excitations were observed to be quite similar. 

%
%
\begin{figure}
\centering
\vspace*{-15mm}

$\begin{array}{l}\mbox{Maximal} \\ \mbox{charge}  
\\ \mbox{transfer} \\[50mm] \textcolor{white}{...} \\
                                    \end{array}$ \;   
\includegraphics[width=0.36\textwidth]{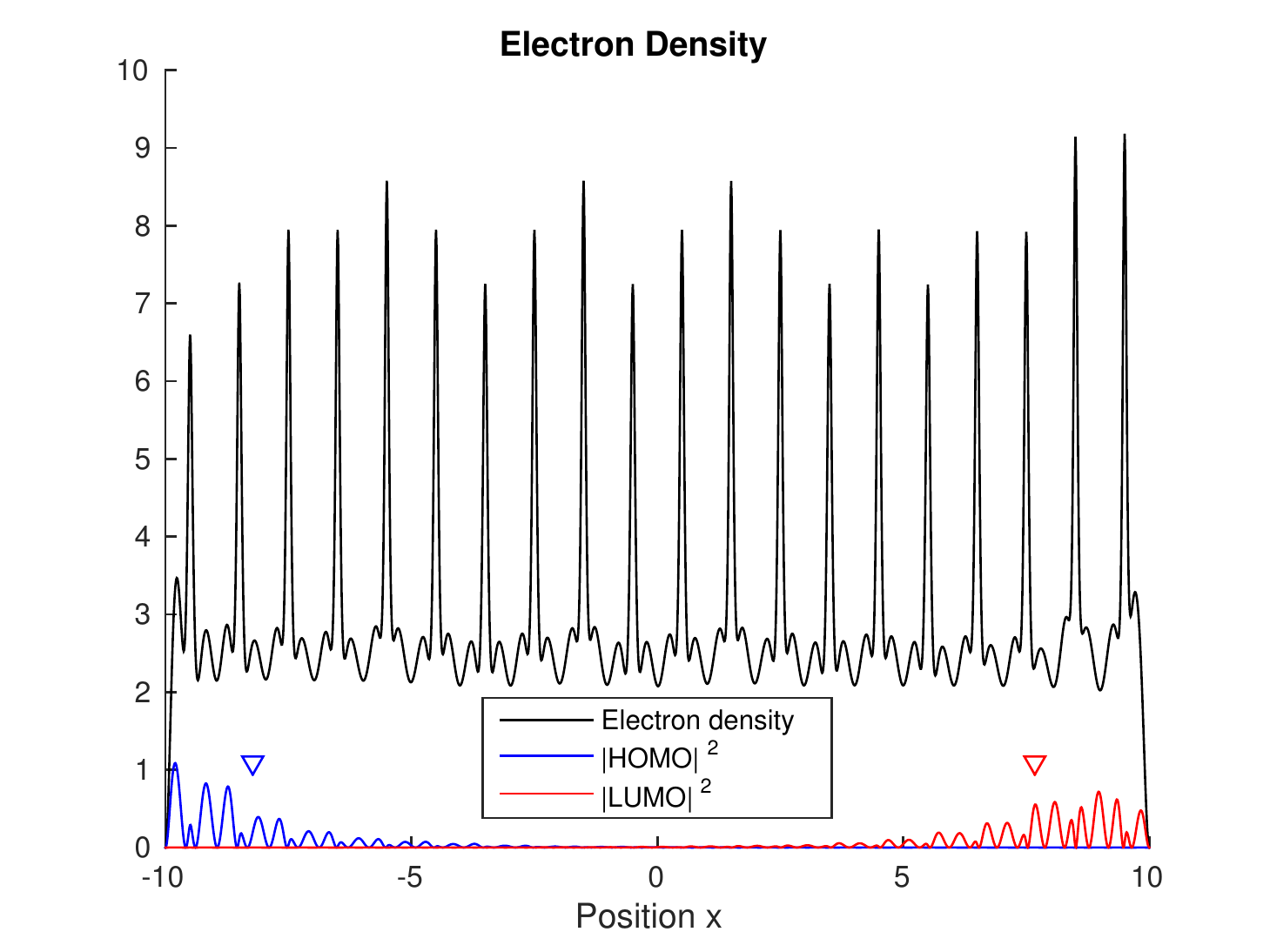} \;
\includegraphics[width=0.36\textwidth]{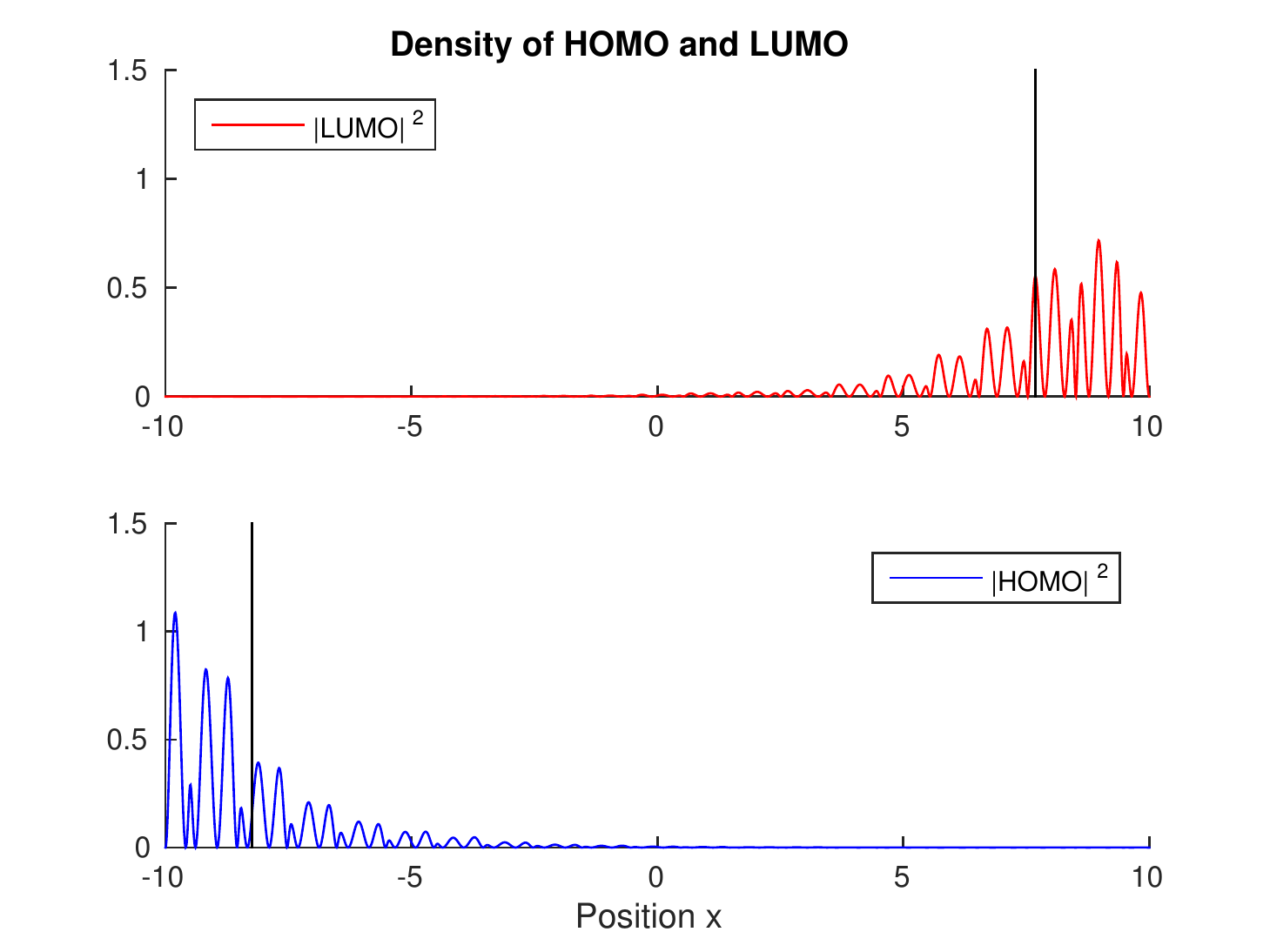}
\vspace*{-28mm} 

$\begin{array}{l}\mbox{Minimal} \\ \mbox{overlap}  
\\ \textcolor{white}{...} \\[50mm] \textcolor{white}{...} \\
                                    \end{array}$ \; \, 
\includegraphics[width=0.36\textwidth]{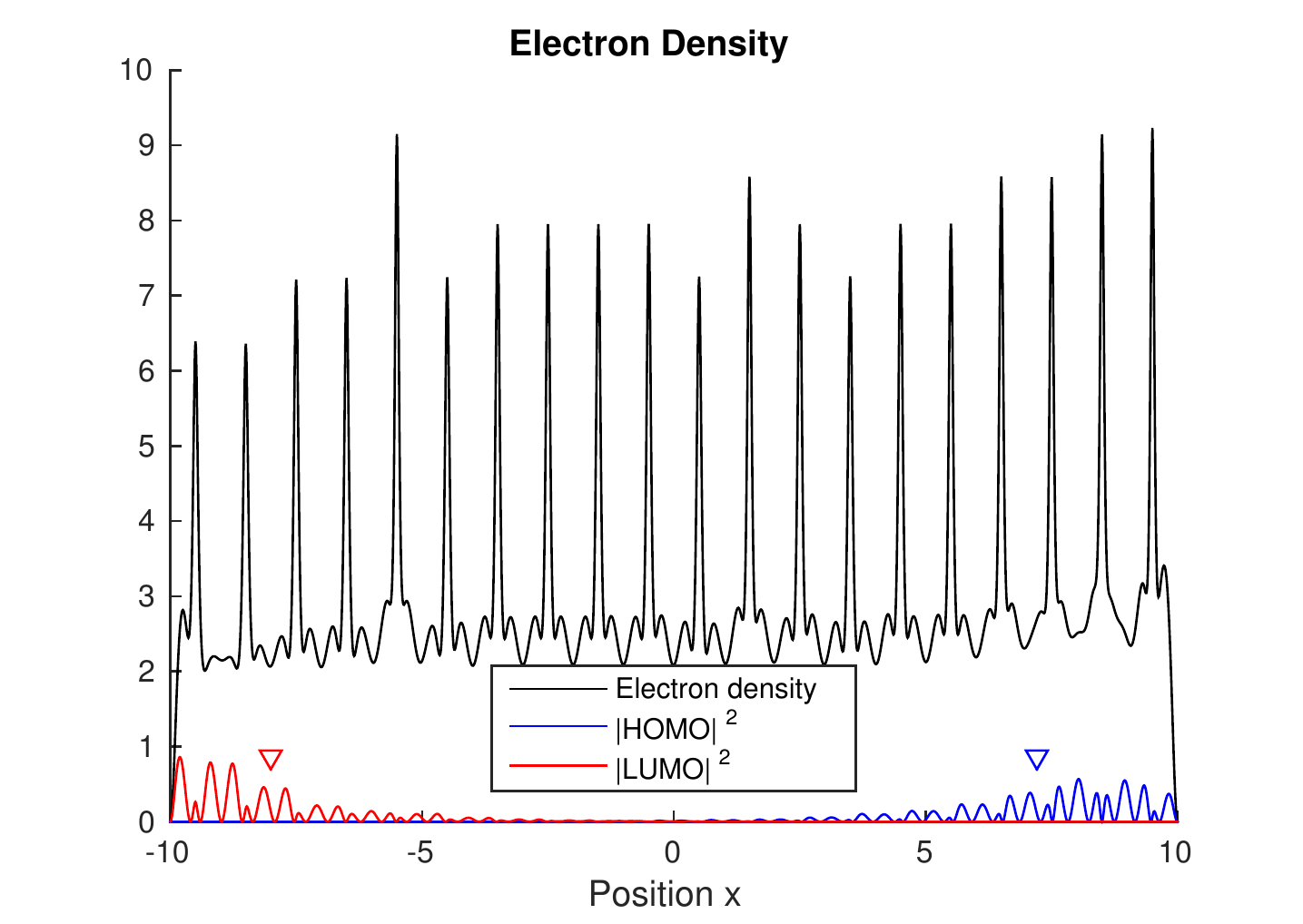} \;
\includegraphics[width=0.36\textwidth]{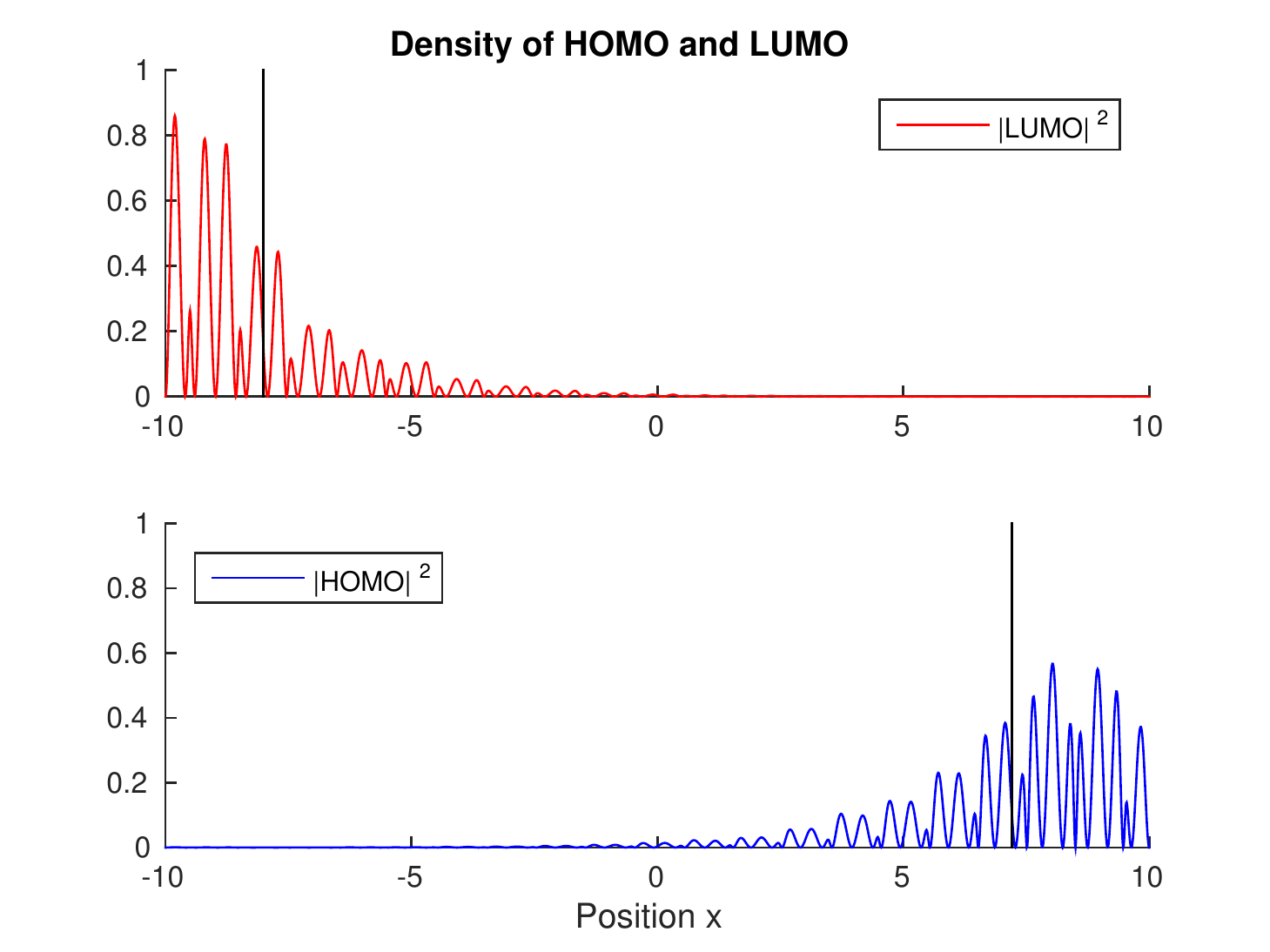}
\vspace*{-28mm} 

$\begin{array}{l}\mbox{Maximal} \\ \mbox{lifetime}  
\\ \textcolor{white}{...} \\[50mm] \textcolor{white}{...} \\
                                    \end{array}$ \,   
\includegraphics[width=0.36\textwidth]{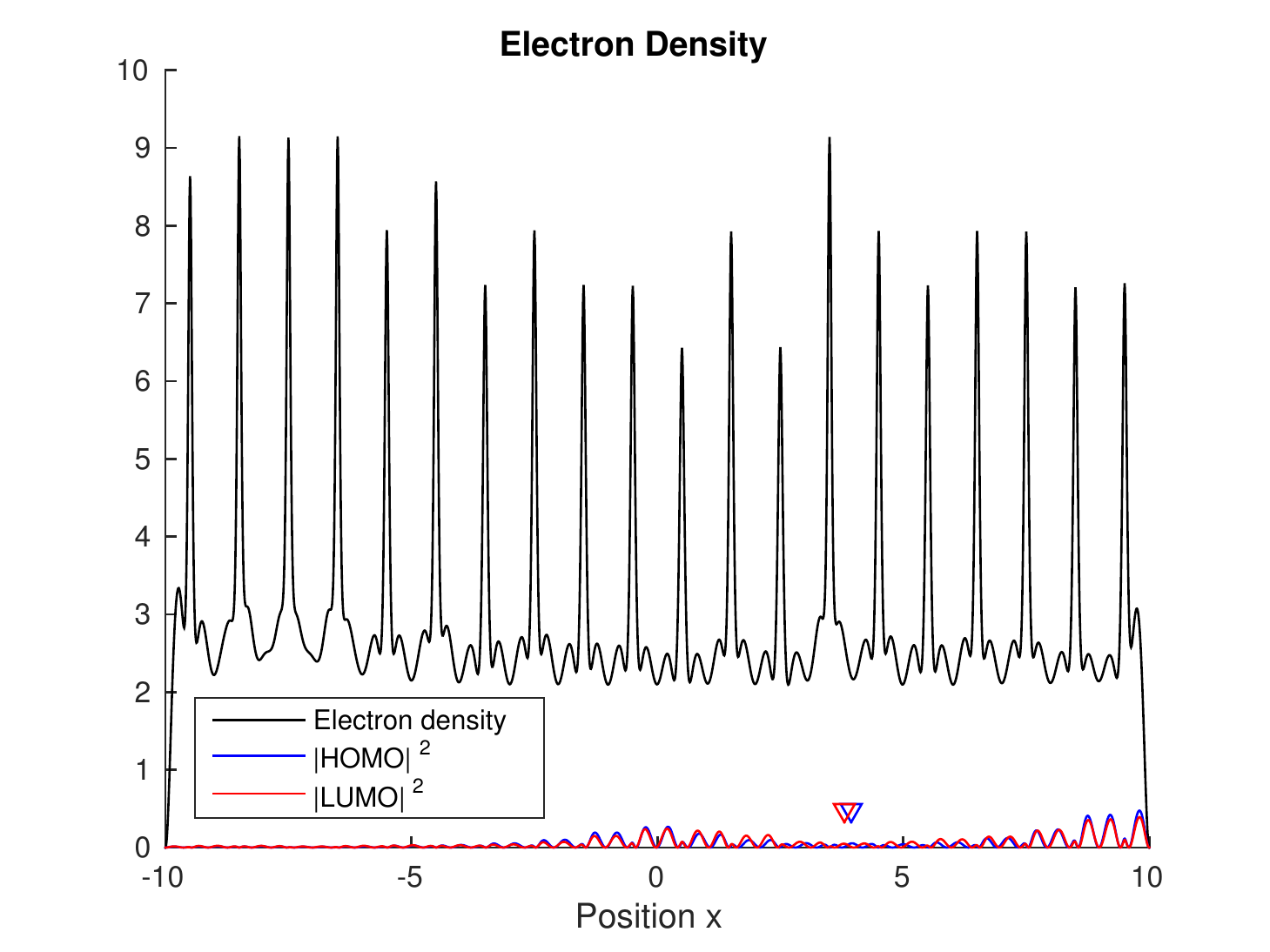}
\includegraphics[width=0.36\textwidth]{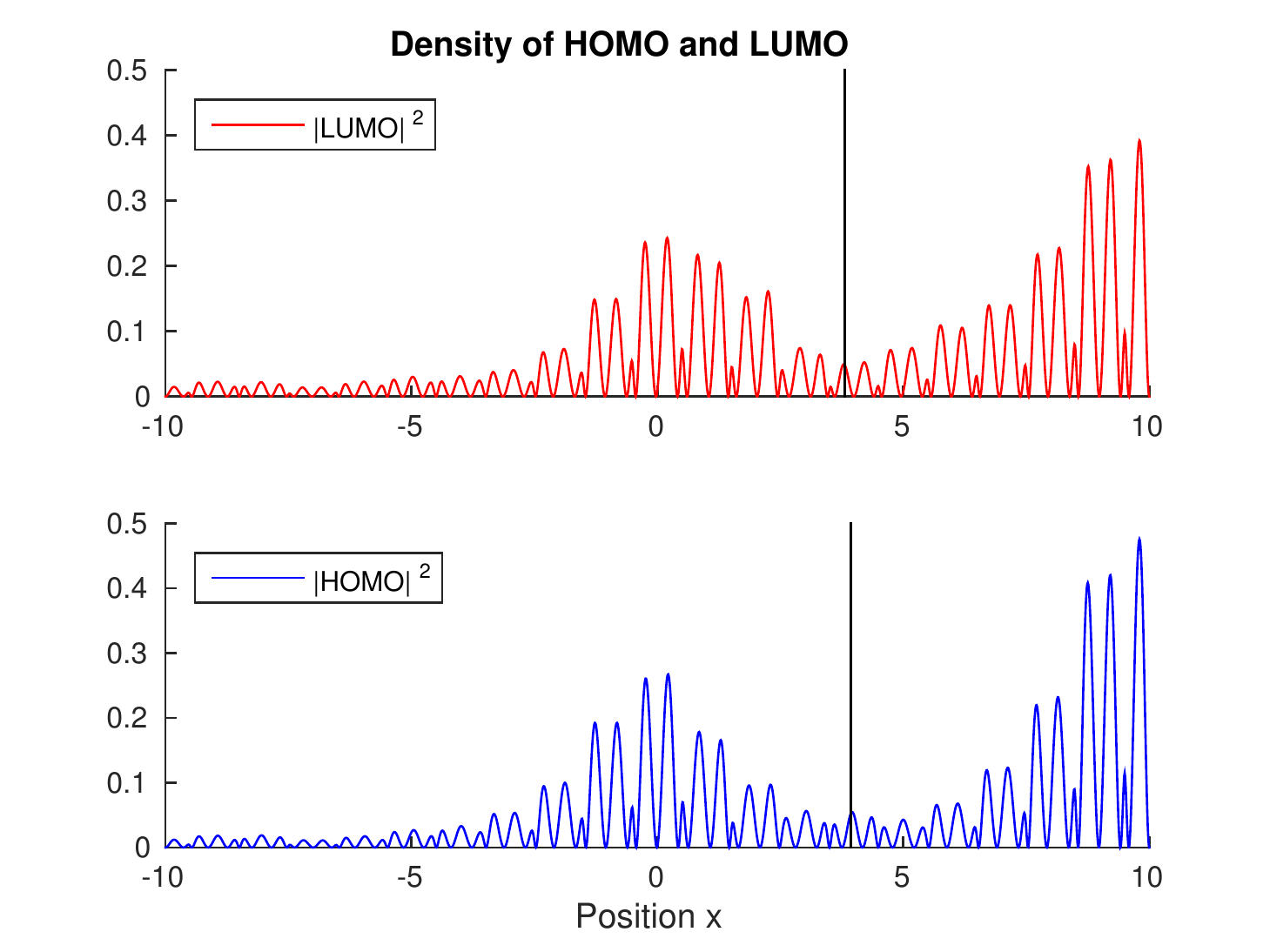}
\vspace*{-28mm} 

$\begin{array}{l}\mbox{Prescribed} \\ \mbox{bandgap}  
\\ \textcolor{white}{...} \\[50mm] \textcolor{white}{...} \\
                                    \end{array}$ \hspace*{-4pt}   
\includegraphics[width=0.36\textwidth]{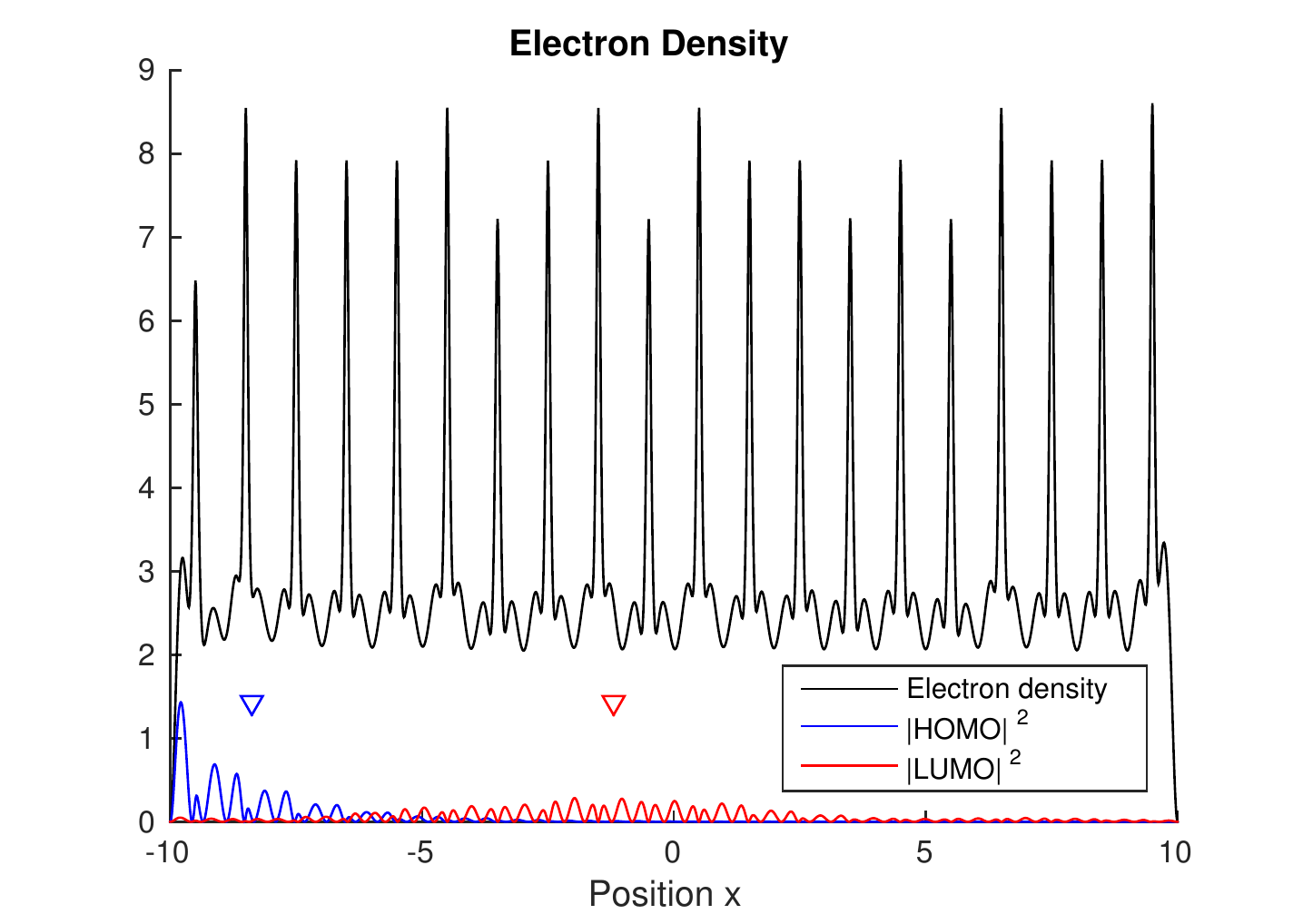}
\includegraphics[width=0.36\textwidth]{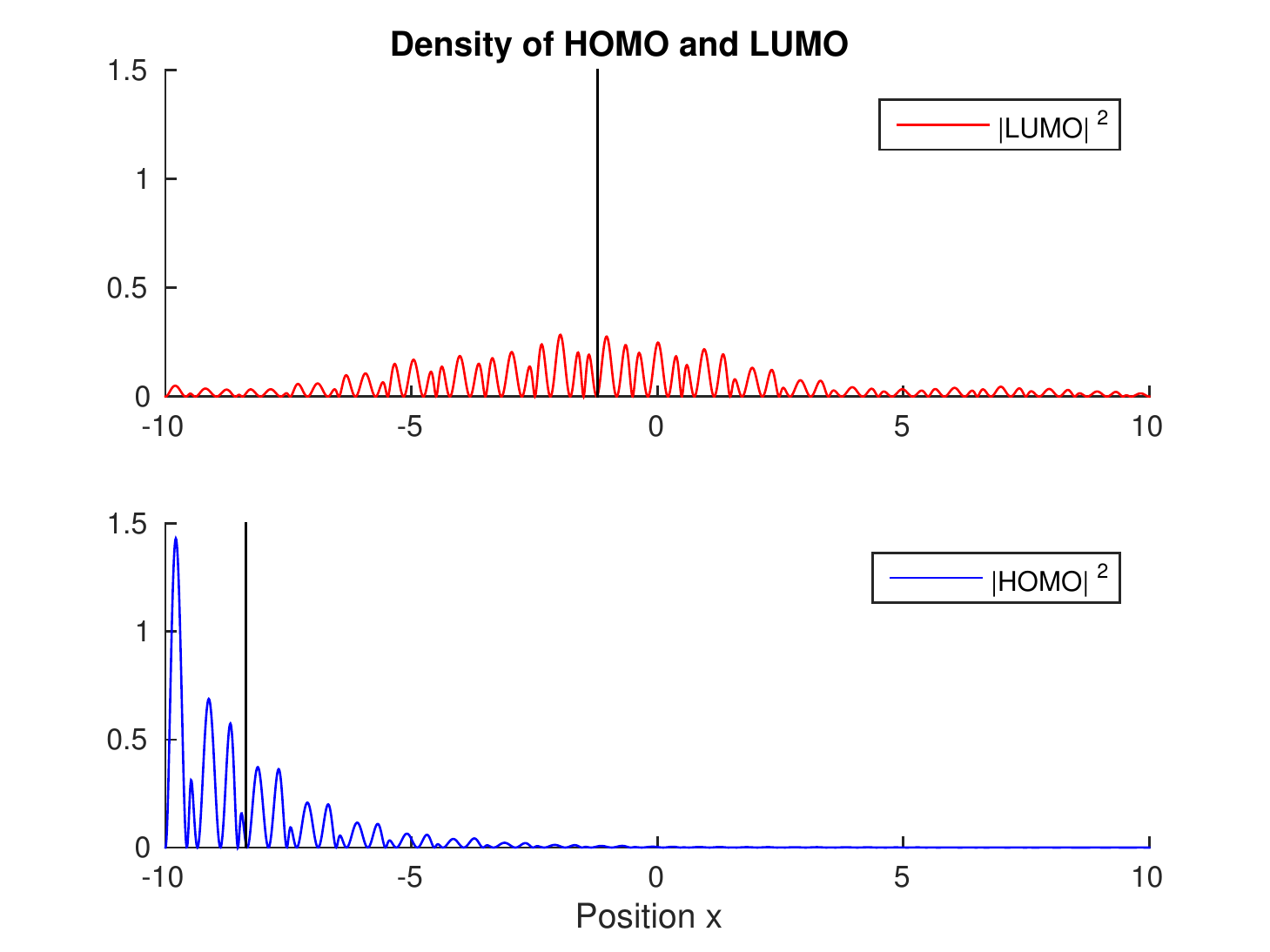}
\vspace*{-33mm}

\caption{{\small HOMO-LUMO excitations obtained by optimization of various goal functionals. For the charge transfer and overlap functionals (first two rows), HOMO and LUMO are well separated on different sides of the system and close to zero on the opposite half of the domain. By contrast, lifetime maximization (third row) leads to nonuniform but nearly identical shapes of HOMO and LUMO; for a simple physical explanation see the text. Finally (bottom row), favouring a prescibed bandgap puts no obvious bias on HOMO and LUMO location and shape, but the prescribed gap of 3 a.u. was reached to high precision (3.0019 a.u.). The optimal doping profiles can be read off from the heights of the electronic density peaks, and exhibit irregular fluctuations.}}
\label{figopt}
\end{figure}

{\bf Pure Carbon Chain.}
The chain of 20 carbon atoms was used as the initial configuration for all our optimization procedures. Figure \ref{figcarbon} depicts the total electron density of the ground state and the density of HOMO and LUMO. On the left, the densities of HOMO and LUMO are not scaled, so as to faithfully indicate their contribution to the overall density. HOMO and LUMO are seen to be delocalized, plane-wave-like, and symmetric with respect to the midpoint of the chain; in particular the net charge transfer \eqref{charge1D} of the HOMO-LUMO-excitation is zero. 
The KS energy levels (see the right panel) exhibit a typical band structure, with a pronounced bandgap between HOMO and LUMO. Also, we encounter 20 very low, near-identical energy levels due to the core states. 
\vspace*{2mm}

\noindent
{\bf Optimal excitations.}
Figure \ref{figopt} shows the excitations achieved by optimization of the four goal functionals \eqref{charge1D}, \eqref{over}, \eqref{life}, and \eqref{designgap}. From the height of the electronic density peaks on the left one can also read off the underlying nuclear configuration, e.g., in case of optimal charge transfer, $(4, 5, 6, 6, 7, 6, 5, 6, 7, 5, 6, 7, 6, 5, 6, 5, 6, 6, 8, 8)$. We conclude from the Figure that our optimal control approach is indeed capable of producing doping profiles whose excitations have the desired features, such as a large charge transfer (first row in Figure \ref{figopt}). 

Next we discuss the nonuniform but nearly identical shapes of HOMO and LUMO in case of lifetime maximization ($3^{rd}$ row in Figure \ref{figopt}), i.e. minimization of \eqref{life}. A simple physical explanation can be given as follows. The nearly identical overall shape leads to an almost vanishing difference in Hartree potential between ground and excited state. Hence the Kohn-Sham Hamiltonian \eqref{ham1D} is nearly identical in the ground and excited state. But the excited state is invariant (up to phase factors) under time evolution with the ground state KS Hamiltonian, and thus almost invariant under time evolution with the excited-state Hamiltonian.

\begin{figure}
\vspace*{-20pt}
\centering
\includegraphics[width=0.5\textwidth]{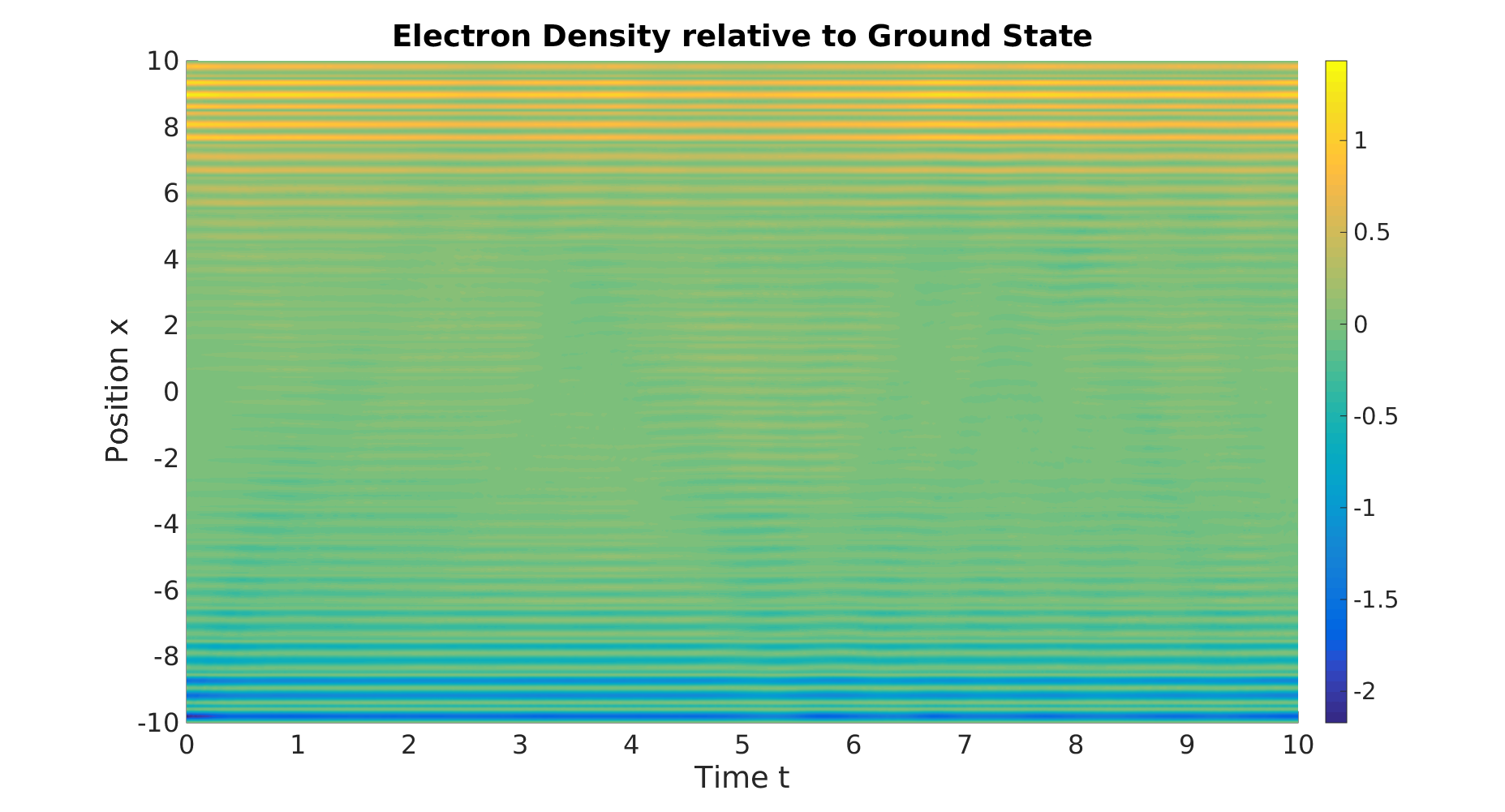} \; \includegraphics[width=0.37\textwidth]{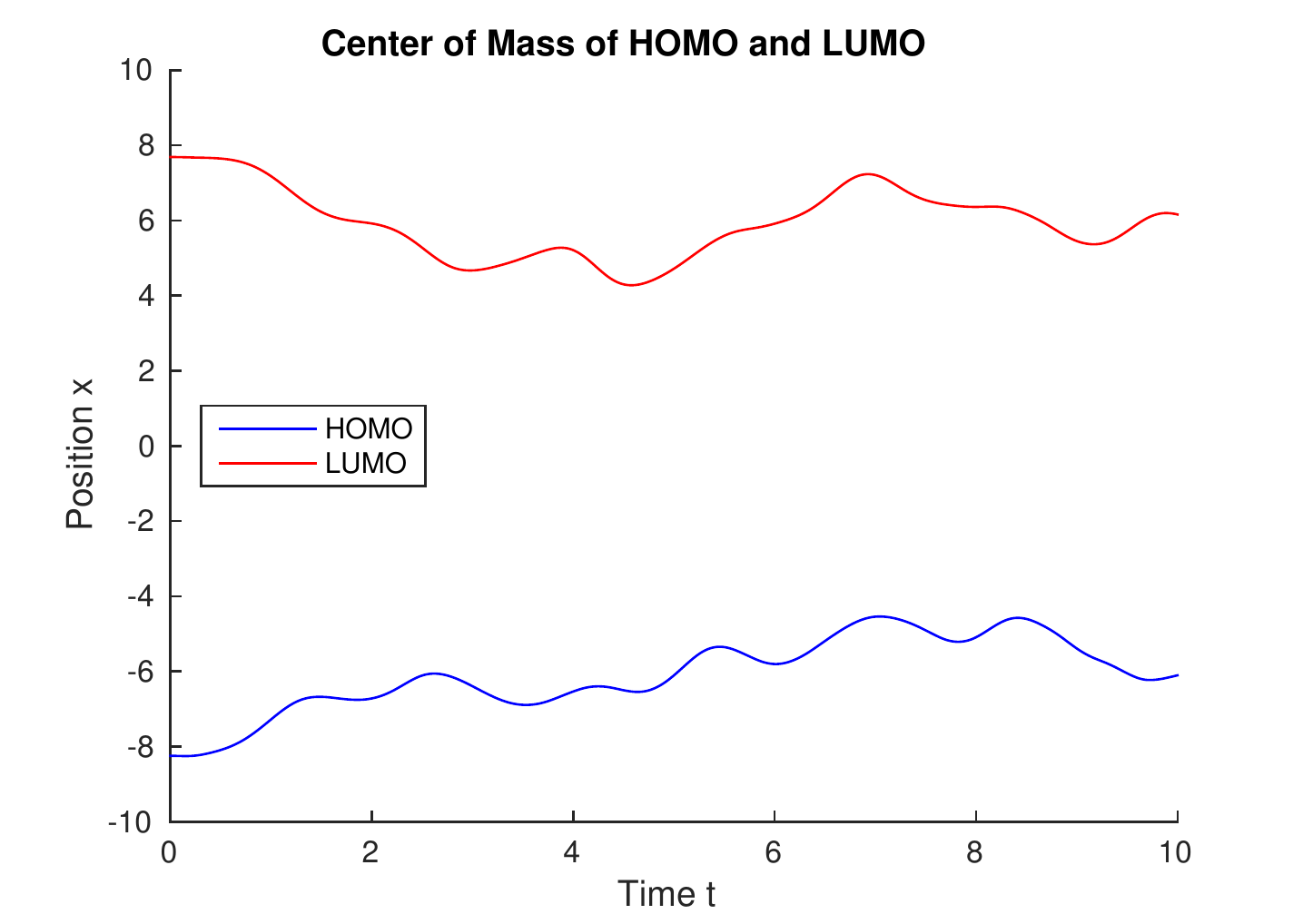} 
\caption{{\small Time evolution of the excitation with maximal charge transfer under TDDFT. Left: electron density relative to the ground state. Right: Center of mass of $\mathrm{HOMO}$ and $\mathrm{LUMO}$.}}
\label{fig:evol}
\end{figure}

Finally let us comment on maximization of the bandgap functional \eqref{gap}. 
Our optimization algorithm didn't find any larger bandgap than that for the pure carbon chain 
($\eps_L - \eps_H = 4.93$); configurations with nearly as high bandgaps found by the algorithm had a large interior pure carbon region and some heteroatoms near the boundary (e.g., $\eps_L-\eps_H=4.88$ for the doping profile $75748566666666577476$). Unfortunately we cannot offer a theoretical explanation for this interesting observation.  
\\[2mm]
{\bf Stability of excitations under time evolution.}
In photovoltaic devices, it is important that excitations persist for a sufficiently long time for harvesting; i.e. one requires a low electron--hole recombination rate. As a minimal check on this we evolved the excitations under TDDFT, eq. \eqref{TDDFT}, with initial conditions given by the new orbitals after excitation. This nonlinear evolution equation (recall that the Hartree potential in \eqref{TDDFT} comes from the time-dependent density) takes the strong electrostatic electron--hole forces fully into account. Figure \ref{fig:evol} depicts, in case of the excitation with maximal charge transfer, the time evolution of the density relative to the ground state and of the center of mass of  $\mathrm{HOMO}$ and $\mathrm{LUMO}$. The results indicate that the electron-hole pair persists during the simulated period of time.

\begin{figure}
\centering
\vspace*{-0pt} 
\hspace*{-.1475\textwidth}
\includegraphics[width=0.37\textwidth]{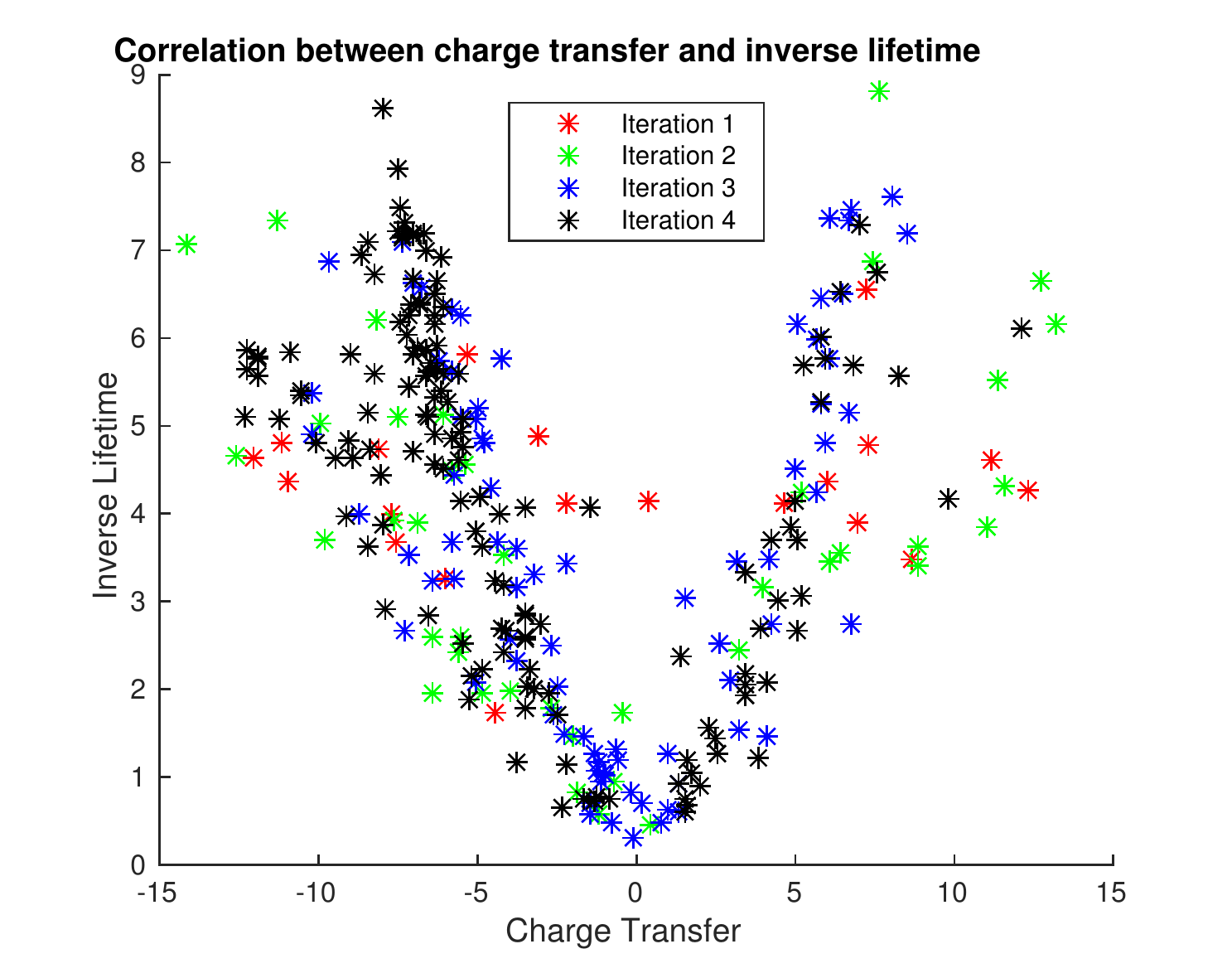}
\hspace*{-.045\textwidth}
\includegraphics[width=0.37\textwidth]{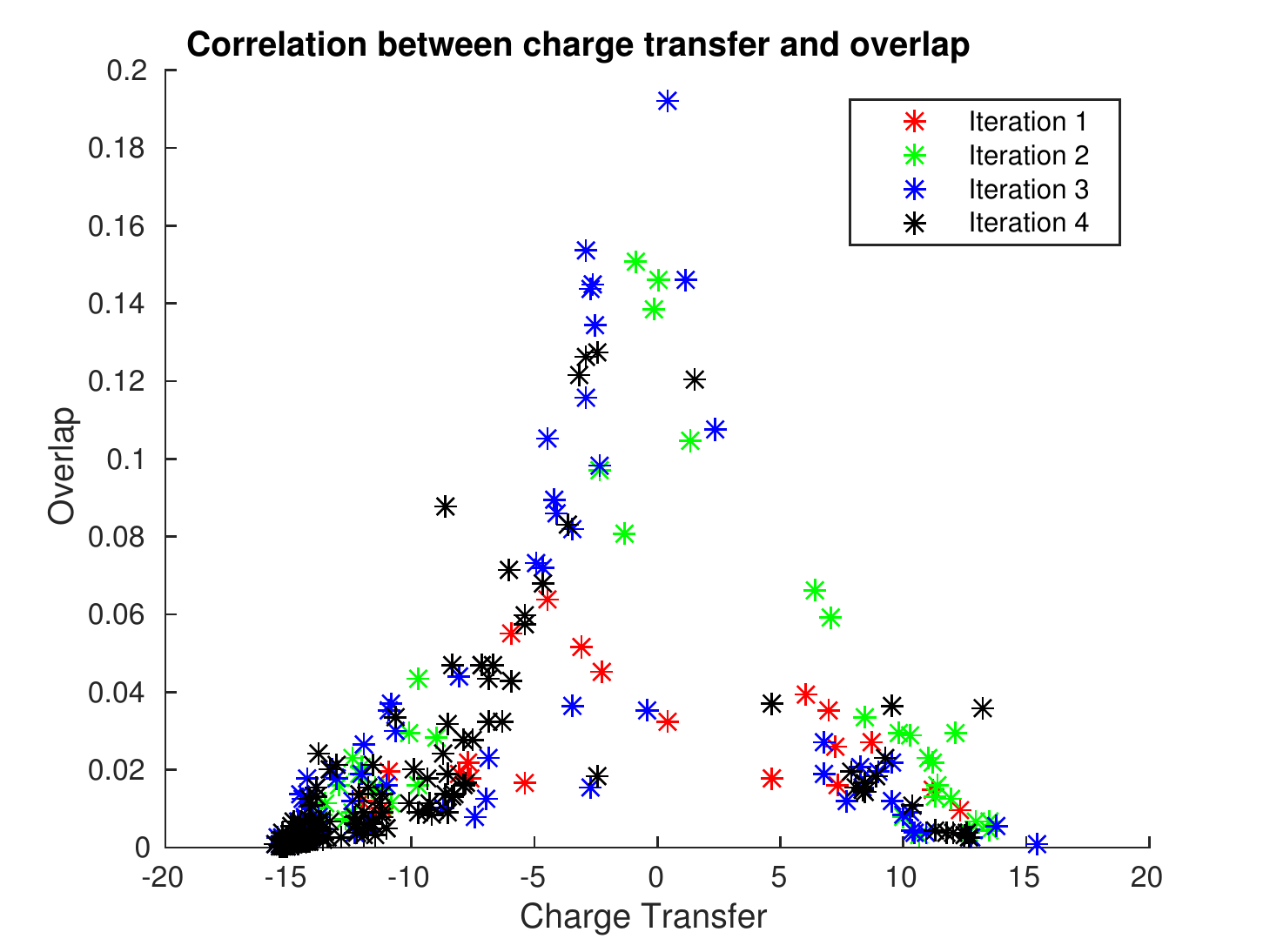}
\hspace*{-.045\textwidth}
\includegraphics[width=0.37\textwidth]{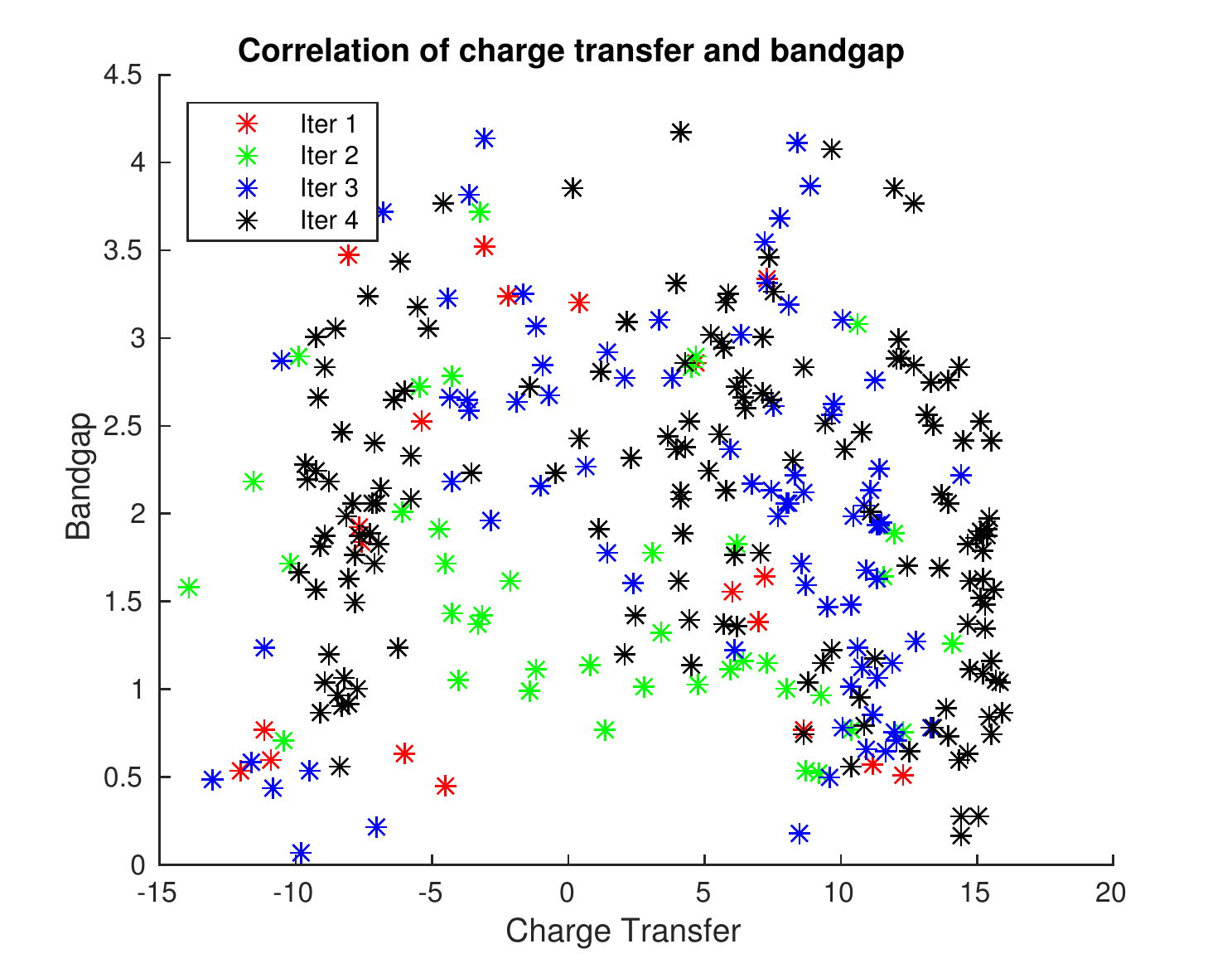}
\hspace*{-.15\textwidth}
\caption{{\small Scatter plots of charge transfer versus inverse lifetime (left), overlap (middle), and bandgap (right), for all nuclear configurations generated during the optimization of 
the quantity on the vertical axis. Coordinates represent numerical values of the functionals \eqref{charge}, \eqref{life}, \eqref{over}, \eqref{gap}.}}
\label{fig:corr}
\end{figure}

\noindent
{\bf Correlations between different excitation properties.} 
When optimizing the nuclear configuration with respect to a property different from charge transfer, we may generate a scatter plot as shown in Figure \ref{fig:corr} where each point corresponds to a configuration generated during our optimization algorithm. Note the strong parabola-like correlation between charge transfer and inverse lifetime (left panel), especially for configurations generated in the last two iteration steps. This shows that a large charge transfer corresponds to a large value of the inverse lifetime functional. We interpret this not as a physical effect but merely as an indication of the limitations of the HOMO-LUMO model: as already mentioned above, a large charge transfer causes a large difference in Hartree potential between ground and excited state, and hence a larger error of the assumption underlying the HOMO-LUMO model that relaxation effects can be neglected. By contrast, we believe that the other two plots capture physical phenomena. First, note the simple inverse correlation between charge transfer and overlap (middle panel). Second,  bandgap and charge transfer appear to be comparatively uncorrelated (right panel), suggesting that these two properties could be controlled simultaneously.  

\vspace*{2mm}

\noindent
{\bf Acknowledgements.} This project was supported by DFG through IGDK 1754 {\it Optimization and Numerical Analysis for Partial Differential Equations with Nonsmooth Structures}. We thank Lin Lin and Jianfeng Lu for helpful discussions. 


\begin{small}

\end{small}


\begin{thebibliography}{99}

\bibitem[AC09]{CancesEtAl} {\sc A. Anantharaman and E. Canc\`{e}s.} {\em Existence of minimizers for Kohn-Sham models in quantum chemistry.} Ann. I. H. Poincar\'{e} 26(6) (2009), 2425–2455.

\bibitem[BSCA03]{BSCA03} {\sc S. Bednarek, B. Szafran, T. Chwiej, and J. Adamowski}, {\em Effective interaction for charge carriers confined in quasi-one-dimensional nanostructures}, Phys. Rev. B 68 (2003), 045328 1--9.

\bibitem[CF15]{CF15} {\sc H. Chen, and G. Friesecke}, {\em Pair densities in density functional theory}, Multiscale Model. Simul. 13(4) (2015), 1259--1289.

\bibitem[FG18]{FG18}{\sc G. Friesecke and B. Graswald}, A proof that the minimizing orbitals of 
some common Kohn-Sham 
energy functionals are lowest eigenstates,
arXiv 2018, to appear

\bibitem[HLYDY16]{HLYDY16} {\sc W. Hu, L. Lin, Ch. Yang, J. Dai, J. Yang.} {\em
Edge-Modified Phosphorene Nanoflake Heterojunctions as Highly
Efficient Solar Cells}, Nano Lett. 16 (2016), 1675−1682

\bibitem[HK64]{HK64} {\sc P. Hohenberg, W. Kohn.} {\em Inhomogeneous electron gas}, { Phys. Rev. B} 136 (1964), 864-871

\bibitem[JN12]{JN12}{\sc R.A.J. Janssen and J. Nelson}, Factors Limiting Device Efficiency in
Organic Photovoltaics, Adv. Mater. 25 (2012), 1847−1858

\bibitem[LK88]{LK88} {\sc Y.-S. Lee and M. Kertesz.} {\em The effect of heteroatomic substitutions on the band gap of polyacetylene and polyparaphenylene derivatives}, J. Chem. Phys. 88 (1988),  2609-2616

\bibitem[KOBH13]{KOBH13}{\sc I.Y. Kanal, S.G. Owens, J.S. Bechtel, and G.R. Hutchison}, J. Phys. Chem. Lett. 4 (2013), 1613−1623


\bibitem[KS65]{KS65} {\sc W. Kohn,  L. J. Sham.} {\em Self-consistent equations including exchange and correlation effects}. { Phys. Rev. A} 140 (1965), pp. 1133-1138 

\bibitem[MW14]{MW14} {\sc J. Ma and L.-W. Wang.} {\em Nanoscale Charge Localization Induced by Random Orientations of Organic Molecules in Hybrid Perovskite CH3NH3PbI3}. { Nano Lett.} 15 (2015), 248−253

\bibitem[P17]{P17} {\sc J.P. Perdew, W.Yang, K.Burke, Z.Yang, E.K.U. Gross, M. Scheffler, G.E. Scuseria, T.M. Henderson, I.Y. Zhang, A. Ruzsinszky, H. Peng, J. Sun, E. Trushin, and A. G\"orling.} {\em Understanding band gaps of solids in generalized Kohn–Sham theory}. PNAS 114 (11) 2801-2806, 2017

\bibitem[PY89]{PY89} R.G. Parr and W. Yang, {\em Density functional theory of atoms and molecules},  
Oxford University Press, 1989.

\bibitem[ZLL11]{ZLL11} {\sc X. Zhang, Z. Li, and G. Lu.} {\em 
First-principles simulations of exciton diffusion in organic semiconductors}. { Phys. Rev. B} 84 (2011), 235208 

\end{thebibliography}
\end{document}